\documentclass[12pt]{article}
\usepackage{amssymb,amsmath,natbib}
\RequirePackage{natbib}
\usepackage{graphicx}
\usepackage{amssymb}
\usepackage{amsthm}
\usepackage{multirow}
\usepackage{array,booktabs}
\usepackage{float}
\usepackage{paralist} 
\usepackage{mathrsfs}
\usepackage{algorithm} 
\usepackage{algorithmic} 
\usepackage[toc,page,title,titletoc,header]{appendix}
\usepackage{authblk}

\newcolumntype{P}[1]{>{\centering\arraybackslash}p{#1}}
\newtheorem{theorem}{Theorem}
\newtheorem{lemma}{Lemma}

\newcommand{\bsb}{\boldsymbol}
\newcommand{\bsbb}{\boldsymbol{b}}
\newcommand{\bsbX}{\boldsymbol{X}}
\newcommand{\bsbZ}{\boldsymbol{Z}}
\newcommand{\bsbz}{\boldsymbol{z}}
\newcommand{\bsbPhi}{\boldsymbol{\Phi}}
\newcommand{\bsbXi}{\boldsymbol{\Xi}}
\newcommand{\bsbx}{\boldsymbol{x}}
\newcommand{\bsbOmega}{\boldsymbol{\Omega}}
\newcommand{\mdiag}{\mbox{diag}}
\newcommand{\bsbLambda}{\boldsymbol{\Lambda}}
\newcommand{\bsblambda}{\boldsymbol{\lambda}}
\newcommand{\bsbDelta}{\boldsymbol{\Delta}}
\newcommand{\bm}{\boldsymbol}
\newcommand{\msgn}{\mbox{sgn}}
\newcommand{\bsbeps}{\boldsymbol{\varepsilon}}
\newcommand{\Proj}{{\boldsymbol{\mathcal P}}}
\newcommand{\EP}{\,\mathbb{P}}
\newcommand{\EE}{\,\mathbb{E}}
\DeclareMathOperator{\vect}{\mbox{vec}\,}
\providecommand{\keywords}[1]{\textbf{Keywords: } #1}

\usepackage[normalem]{ulem}

\begin{document}

\title{Group Regularized Estimation under Structural Hierarchy}
\author{Yiyuan She,    Zhifeng Wang and He Jiang}
\affil{Department of Statistics, Florida State University}
\date{}
\maketitle


\begin{abstract}
 Variable selection for models including interactions between explanatory variables often needs to obey certain hierarchical constraints. The weak or strong structural hierarchy requires that the existence of an interaction term implies at least one or both associated main effects to be present in the model.
Lately, this problem has attracted a lot of attention, but existing computational algorithms converge slow even with a moderate number of predictors. Moreover, in contrast to the rich literature on ordinary variable selection, there is a lack of statistical theory to show reasonably low error rates of hierarchical variable selection.
 This work investigates a new class of estimators that make use of multiple group    penalties to capture structural parsimony.
We give  the minimax lower bounds for strong and weak hierarchical variable selection and show that the proposed estimators enjoy sharp rate oracle inequalities. A general-purpose algorithm is developed with guaranteed convergence and global optimality. Simulations and real data experiments  demonstrate the efficiency and efficacy of the proposed approach.

\end{abstract}



\section{Introduction}
\label{sec:intro}

In statistical applications, it is often noticed that an additive model including main effects only is inadequate. Including some higher-order terms, such as interactions, in particular,  are often of great help in prediction and modeling. Sometimes, interactions may be of independent interest; one example is the moderation analysis in behavioral sciences \citep{cohen2013applied}. In this paper, we focus on the full quadratic model with all two-term interactions taken into account.

Let $\bsbX =[\bsb{x}_1,\bsb{x}_2,\ldots,\bsb{x}_p] \in \mathbb{R}^{n \times p}$ be the (raw) predictor matrix and $\bsb{y} \in \mathbb{R}^n$ be the response vector. We assume the following nonlinear additive regression model
\begin{equation}
\bsb{y}=b_0^* \bsb{1}+\sum_{1\leq j \leq p} b_j^* \bsb{x}_j+\sum_{1 \leq j, k \leq p}\phi_{jk}^* \bsb{x}_j\odot\bsb{x}_k+\bsb{\varepsilon}, \label{model}
\end{equation}
where $\bsb{\varepsilon} \sim N(\bsb{0}, \sigma^2\bsb{I})$ and $\odot$ denotes the Hadamard product. The number of predictors is then $p+{p \choose 2}$, posing a challenge in variable selection even when $p$ is moderate.
Moreover, in this scenario, statisticians are often interested in obtaining a model satisfying certain logical relations, such as the {\textit{structural hierarchy}} discussed in   \cite{Nelder}, \cite{SH}, and \cite{WH}.  Hierarchy is a natural requirement in  gene regulatory network studies \citep{davidson2006gene}, banded covariance matrix estimation \citep{bien2014convex} and lagged variable selection in time series.
Hierarchical variable selection leads to reduced number of variables in measurement, referred to as \emph{practical sparsity} \citep{HL}.  For instance, a model consisting of $\bsbx_1$, $\bsbx_2$, and $\bsbx_1\bsbx_2$  may be more parsimonious to practitioners than a model involving $\bsbx_1$, $\bsbx_2$, and $\bsbx_3\bsbx_4$.
In our setting, there are two types of hierarchy \citep{Chipman, HL}: strong hierarchy (\textbf{SH}) and {weak hierarchy} (\textbf{WH}). Let $\phi_{jk}$ be the coefficient of  $\bm{x}_j\odot\bm{x}_k$ and $\bm{x}_k\odot\bm{x}_j$. SH means that if an interaction term exists in the model, then both of its associated main effects must be present, i.e., $\phi_{jk} \neq 0  \rightarrow b_j \neq 0 \ \text{and} \ b_k \neq 0$, while WH requires that the inclusion of an interaction implies   at least one of its associated main effects to be added into the model, i.e., $\phi_{jk} \neq 0  \rightarrow b_j \neq 0 \ \text{or} \ b_k \neq 0$.
We will show in Section \ref{sec:greshframework} that WH is relatively easy to realize compared with SH. SH is invariant to linear transformations of  predictors \citep{peixoto1990property} and is the primary concern in this work.

It is a nontrivial task to maintain hierarchy in model selection  using conventional approaches. LASSO \citep{LASSO} may violate SH and WH as well. We refer to \cite{Nelder}, \cite{Peixoto},  \cite{bickel2010hierarchical}, \cite{wu2010screen}, and \cite{hao2014interaction} for some well-developed multi-step procedures which, however, might be ad hoc and greedy.
This paper focuses on \emph{regularization-based} approaches. The past works in this direction  include SHIM \citep{CHOI}, VANISH \citep{VANISH} and HL \citep{HL}. SHIM reparametrizes $\phi_{jk}$ as $\rho_{jk}b_jb_k$ and enforces sparsity in both $\bsbb = [b_j]$ and $\bsb{\rho} = [\rho_{jk}]$. The formulation is motivating, and we could also use     $\phi_{jk} = \rho_{jk}(b^2_j + b^2_k)$ for WH.
 However, the corresponding optimization problem is nonconvex and the computational algorithm of SHIM is quite slow in large-scale problems. VANISH is one of the main motivations of our work and will be discussed in detail in Section \ref{sec:greshframework}. {HL is a recent breakthrough in hierarchical variable selection}. One of its key ideas is to {enforce} a magnitude constraint on the coefficients, $\|\bsb{\phi}_j\|_1 \leq |b_j|$, to make hierarchy naturally hold. Here, $\bsb{\phi}_j$ is a vector of coefficients of the predictors $\bsbx_j \odot \bsbx_k$, $1 \leq k \leq p$. To handle the nonconvex constraint, \cite{HL} rephrased it with   the pesudo-positive and pesudo-negative parts $b_j^+$, $b_j^-$ of $b_j$  but dropped all zero-product constraints $b_j^+ b_j^-=0, 1\le j\le p$. The quality of such a convex relaxation seems to have no theoretical justification in the literature. In our experience, HL has excellent performance when the main effects are strong and $p$ is not very large. But it can miss some interaction effects  and   become computationally prohibitive on  large datasets. For example, when $p=1000$, HL can take days to obtain a 20-point solution path.

In this work, we propose and study group regularized estimation under structural hierarchy (\textbf{GRESH}). In theory, we are able to establish non-asymptotic oracle inequalities to show the error rates of the proposed estimators are   minimax optimal up to some logarithm factors. We come up with a new recipe to conquer the theoretical difficulties when analyzing overlapping regularization terms in pursuing structural parsimony.  Moreover, we develop a  computational algorithm which guarantees the convergences of iterates and function values; it is  not only efficient but also simple to implement.

The rest of the paper is organized as follows. Some notation and symbols are introduced in Section \ref{sec:notation}. Section \ref{sec:greshframework} presents the general framework of GRESH. A fast computational algorithm with theoretical support is given in Section \ref{sec:comp}. Section \ref{sec:oracle} builds oracle inequalities for GRESH, and Section \ref{sec:minimax} shows the minimax optimal rates.
In Section \ref{sec:exp}, simulation studies and real data analysis are conducted to show the prediction accuracy and computational efficiency of the proposed approach. All technical proofs  are given in the Appendices.

\section{Notation}
\label{sec:notation}

We introduce some convenient notation and symbols to be used throughout the paper. First, for any matrix $\bsb{A}=\left[\bsb{a}_1,\ldots,\bsb{a}_p\right]^{\top} \in \mathbb{R}^{p \times n}$, define its $(2,1)$-norm, $(2,
\infty)$-norm and $\ell_1$-norm as $\|\bsb{A}\|_{2,1} = \sum_{i=1}^p \|\bsb{a}_i\|_2$, $\|\bsb{A}\|_{2,\infty}=\max_{1 \leq i \leq p}\|\bsb{a}_i\|_2$ and $\|\bsb{A}\|_1= \|\vect(\bsb{A})\|_1$, respectively, where $\vect$ is the standard vectorization operator. The spectral norm and Frobenius norm of $\bsb{A}$ are denoted by $\|\bsb{A}\|_2$ and $\|\bsb{A}\|_F$, respectively. For any $p$-dimensional vector $\bsb{a}$ that is divided into $K$ groups with $\bsb{a}_j$ representing the $j$-th subvector, its $(2,1)$-norm is defined by $\|\bsb{a}\|_{2,1}=\sum_{j=1}^{K}\|\bsb{a}_j\|_2$.

The following two operators $\mdiag$ and $\mbox{dg}$ are introduced for notational simplicity. For a square matrix $\bsb{A}:=[a_{ij}]_{n\times n}$, $\mbox{diag}(\bsb{A}):=[a_{11},\ldots,a_{nn}]^{\top}$, and for a vector $\bsb{a}=[a_1,\ldots,a_n]^{\top}$ $\in \mathbb{R}^{n}$, $\mbox{diag}\{\bsb{a}\}$ is defined as an $n \times n$ diagonal matrix with diagonal entries given by $a_1,\ldots,a_n$. Define  $\mbox{dg}(\bsb{A})  :=\mdiag\{\mdiag(\bsb{A})\}=\mdiag\{[a_{11},\ldots,a_{nn}]^{\top}\}$.
We use $\bsb A[\mathcal I, \mathcal J]$ to  denote a submatrix of $\bsb A$ with rows and columns indexed by $\mathcal I$ and $\mathcal J$, respectively.

For any arbitrary $\bsbb \in \mathbb{R}^p$, $\bsbPhi \in \mathbb{R}^{p \times p}$, we define
\begin{equation}
\begin{aligned}
\mathcal{J}^{11}(\bsbb,\bsbPhi):&=\{j \in [p] : b_j \neq 0,  \ \bsb{\phi}_j \neq \bsb{0}\},\\
\mathcal{J}^{10}(\bsbb,\bsbPhi):&=\{j \in [p]:  b_j \neq 0, \ \bsb{\phi}_j = \bsb{0}\},\\
\mathcal{J}^{01}(\bsbb,\bsbPhi):&=\{j \in [p]:  b_j =0, \ \bsb{\phi}_j  \neq \bsb{0}\},\\
\mathcal{J}^{00}(\bsbb,\bsbPhi):&=\{j \in [p]:  b_j =0, \ \bsb{\phi}_j = \bsb{0}\},\\
\mathcal{J}_e(\bsbPhi):&=\{j\in [p^2]:(\vect(\bsbPhi))_j \neq 0\},\\
\mathcal{J}_G(\bsbb,\bsbPhi):&=\{j\in [p] :b^2_j+\|\bsb{\phi}_j\|_2^2 \neq 0\},
\label{Jdefs}
\end{aligned}
\end{equation}
 where $[p]=\{1,\ldots,p\}$, and the coefficient vector $\bsb{\phi}_j$ denotes the $j$-th column of $\bsbPhi$. With $|\cdot|$ standing for set cardinality, we define
 $J^{11}(\bsbb,\bsbPhi):=|\mathcal{J}^{11}(\bsbb,\bsbPhi)|$, $J^{10}(\bsbb,\bsbPhi):=|\mathcal{J}^{10}(\bsbb,\bsbPhi)|$, $J^{01}(\bsbb,\bsbPhi):=|\mathcal{J}^{01}(\bsbb,\bsbPhi)|$, $J^{00}(\bsbb,\bsbPhi):=|\mathcal{J}^{00}(\bsbb,\bsbPhi)|$, $J_e(\bsbPhi):=|\mathcal{J}_e(\bsbPhi)|$ and $J_G(\bsbb,\bsbPhi)$ $:=|\mathcal{J}_G(\bsbb,\bsbPhi)|$.
 Clearly, $J^{11}+J^{10}+J^{01}+J^{00}=p$, and $J^{11}+J^{10}+J^{01}=J_G$. In addition, under SH, $J_G(\bsbb,\bsbPhi)$ equals the number of nonzero elements of $\bsb{b}$. Given the true signal $(\bsbb^*,\bsbPhi^{*})$, the following abbreviated symbols are used: $J^{11*}=J^{11}(\bsbb^*,\bsbPhi^*)$, $J^{10*}=J^{10}(\bsbb^*,\bsbPhi^*)$, $J^{01*}=J^{01}(\bsbb^*,\bsbPhi^*)
$, $J_{e}^* = J_e(\bsbPhi^*)$, and $J^{*}_{G}=J_G(\bsbb^*,\bsbPhi^*)$.

In the paper, we frequently use the concatenated coefficient matrix for convenience
\begin{equation}
 \bsbOmega=[\bsbb, \bsbPhi^{\top}]^{\top},
 \end{equation}
and its $j$-th column is denoted by $\bsbOmega_j=[b_j,\bsb{\phi}_j^{\top}]^{\top}$. Given $\bsbOmega$,   $\bsbOmega_{b}$ and $\bsbOmega_{\Phi}$  stand for   $(\bsbOmega[1,:])^{\top}$ and  $\bsbOmega[2:(p+1),:]$, respectively.

Let $\bsbX = [\bsbx_1,\ldots, \bsbx_p] \in \mathbb{R}^{n \times p}$ be the raw predictor matrix.
Define
\begin{align}
 \bar{\bsbX} &=[\bsbx_1 \odot \bsbx_1,  \ldots, \bsbx_1 \odot \bsbx_p, \ldots, \bsbx_p \odot \bsbx_1, \ldots, \bsbx_p \odot \bsbx_p]\in \mathbb{R}^{n\times p^2}, \label{defXbar}
\\
\breve{\bsbX}&=[\bsbx_1, \bsbx_1 \odot \bsbx_1, \ldots,\bsbx_1 \odot \bsbx_p,  \ldots, \bsbx_p, \bsbx_p \odot \bsbx_1, \ldots, \bsbx_p \odot \bsbx_p] \in \mathbb{R}^{n\times (p^2+p)}. \label{defXbreve}
\end{align}
Then $\bar{\bsbX}$ consists of all interactions, and $\breve{\bsbX}$ includes all $p^2+p$ predictors in the quadratic model.
Given any subset $\mathcal J \subset [p]$, we abbreviate $\bsb X[:,\mathcal J]$ as $\bsb X_{\mathcal J}$.
It is also easy to see that $\mdiag(\bsbX\bsbPhi\bsbX^{\top})=\bar{\bsbX}\vect(\bsbPhi)$ and $\bsbX\bsbb +\mdiag(\bsbX\bsbPhi\bsbX^{\top})= \breve\bsbX \vect(\bsbOmega)$.

For any two real numbers $a$ and $b$, $a \lesssim b$ means that $a \leq b$ holds up to a multiplicative numerical constant. For two equally sized matrices $\bsb{A}=(a_{ij})$ and $\bsb{B}=(b_{ij})$, $\bsb{A} \geq \bsb{B}$ means $a_{ij} \geq b_{ij}$, for $\forall \ i,j$.

\section{Group regularized estimation under structural hierarchy}
\label{sec:greshframework}

For simplicity, we assume for now that there exists no intercept term in the model. Then \eqref{model} can be written as
\begin{equation}
\bsb{y}=\bsbX\bsbb^*+\mbox{diag}(\bsbX\bsbPhi^*\bsbX^{\top})+\bsb{\varepsilon}, \label{model1}
\end{equation}
where  $\bsb{\varepsilon} \sim N(\bsb{0}, \sigma^2\bsb{I})$ with $\sigma^2>0$, $\bsb{y}=[y_1,\ldots,y_n]^{\top} \in\mathbb{R}^n$ is the response vector, and $\bsb{X}=[\bsbx_1,\ldots,\bsbx_p] \in \mathbb{R}^{n \times p}$ is the design matrix consisting of main effects only.

We describe a general framework  for hierarchical variable selection, referred to as group regularized estimation under structural hierarchy or \textbf{GRESH}. GRESH has two different types, depending on which objects to regularize. Denoting the squared error loss by $\ell$, i.e.,
\begin{equation}
\ell(\bsbb, \bsbPhi)=\frac{1}{2}\|\bsb{y}-\bsbX\bsbb-\mbox{diag}(\bsbX\bsbPhi\bsbX^{\top})\|_2^2, \label{loss}
\end{equation}
the first type is given by
\begin{equation}
\begin{gathered}
\textbf{Type-A}: \min_{\bsbOmega=[\bsbb,\bsbPhi^{\top}]^{\top}\in \mathbb{R}^{(p+1)\times p}}
 \ell(\bsbb, \bsbPhi) +\lambda_1\|\bsbPhi\|_1+
\lambda_2\sum_{j=1}^p\|[b_j,z(\bsb{\phi}_j)] \|_q \\
\text{s.t.} \ \bsbPhi=\bsbPhi^{\top} \ (\text{for \textbf{SH} only}),
\end{gathered}   \label{ob1}
\end{equation}
where $\lambda_1,\lambda_2$ are regularization parameters, $1<q\leq +\infty$ and  $z(\bsbx)$ is a function satisfying the property that $z(\bsbx) =\bsb{0}$ implies $\bsbx=\bsb{0}$ for any vector $\bsbx\in \mathbb{R}^p$. For instance, $z$ can take the  $\ell_r$-norm function ($r>0$)
\begin{equation}
z(\bsbx)=\|\bsbx\|_r, \label{r}
\end{equation}
or simply the identity function
\begin{equation}
z(\bsbx)=\bsbx^{\top}.        \label{identity}
\end{equation}
The first $\ell_1$ penalty imposes elementwise sparsity on $\bsbPhi$ and the second group-$\ell_1$ penalty enforces column sparsity in $\bsbOmega$. We argue that with the two penalties and the constraint, \eqref{ob1} can be used for strong hierarchical variable selection. Indeed, the sparsity of $\bsbb$ comes from the second group-penalty alone, i.e., $b_j=0$ implies $\|[b_j,z(\bsb{\phi}_j)] \|_q=0$ (with probability 1) or $z(\bsb{\phi}_j)=0$. By the property of the $z$ function, $\bsb{\phi}_j=\bsb{0}$, and thus $\phi_{jk}=0$. The symmetry condition indicates further that $\phi_{kj}=0$. Hence $(\phi_{jk}+\phi_{kj})/2$, the coefficient for $\bsbx_j \odot \bsbx_k$,  is zero. Consequently, whenever $b_j = 0 $, $\bsbx_j \odot \bsbx_k$ will be removed from the model and     SH is automatically obeyed. We can describe the reasoning as follows
\begin{equation} \label{derivation}
\begin{gathered}
b_j = 0  \Rightarrow  \|[b_j,z(\bsb{\phi}_j)] \|_q = 0  \Rightarrow  z(\bsb{\phi}_j)=0  \Rightarrow  \bsb{\phi}_j = 0  \Rightarrow  \phi_{jk} = 0\\
\Rightarrow  \phi_{kj} = 0  \Rightarrow  \frac{\phi_{jk}+\phi_{kj}}{2} = 0.
\end{gathered}
\end{equation}
Without the symmetry constraint, we can only complete the argument in the first line of \eqref{derivation}, and so SH dose not hold. But interestingly,   WH is guaranteed, because from $b_j = b_k = 0$, we have $(\phi_{jk}+\phi_{kj})/2=0$. Therefore, WH gives a relatively simpler problem.

As pointed out by a reviewer, when the model contains an intercept, centering the response and the $p$ raw predictors does not make it vanish due to the presence of nonlinear terms, and so  $\bsb{1} b_0 + \bsbX\bsbb+\mbox{diag}(\bsbX\bsbPhi\bsbX^{\top})$ should be used to approximate $\bsb{y}$. In the SH scenario, if  at least one $x$-predictor is relevant, substituting $\bsbX'=[\bsb{1}, \bsbX]$ for $\bsbX$ in  \eqref{loss}  suffices.

We focus on convex forms of GRESH in this work. But surely the $\ell_1$ penalty and the  group-$\ell_1$ penalty in \eqref{ob1} can be replaced by their nonconvex alternatives; see, e.g., \cite{Group_She}.

GRESH is related to some methods in the literature. HL makes a special case of \eqref{ob1} because one of its formulations corresponds to $q=\infty$ and $r = 1$, with a single regularization parameter being used. Another instance is given by $q=2$ and $z(\bsbx)=\bsbx^{\top}$:
\begin{equation}
\min_{\bsbOmega=[\bsbb,\bsbPhi^{\top}]^{\top}}\ell(\bsbb, \bsbPhi) +\lambda_1\|\bsbPhi\|_1+\lambda_2\|\bsbOmega^{\top}\|_{2,1} \quad  \text{s.t.} \ \bsbPhi=\bsbPhi^{\top}. \label{typeA_comm}
\end{equation}
\cite{HL} incorrectly described VANISH \citep{VANISH} in this form, without the symmetry condition.  We will focus on \eqref{typeA_comm} in the theoretical and computational studies of Type-A GRESH.

 As a matter of fact, \cite*{VANISH} defined VANISH in a different way, which motivates another type of GRESH
\begin{equation}
\begin{gathered}
\textbf{Type-B:}  \min_{\bsbOmega=[\bsbb,\bsbPhi^{\top}]^{\top }}\ell(\bsbb, \bsbPhi)+ \lambda_1 \sum_{1 \leq j,k \leq p}\|\phi_{jk}\bsbx_j \odot \bsbx_k \|_2 + \\
\lambda_2\sum_{j=1}^{p}\|[b_j\bsbx_j,z(\bsb{\phi}_{j1}\bsbx_j \odot \bsbx_1,\ldots, \bsb{\phi}_{jp}\bsbx_j \odot \bsbx_p)]\|_q, \, \text{ s.t. }  \bsbPhi = \bsbPhi^{\top} (\mbox{\textbf{SH} only}),
\end{gathered} \label{Type B}
\end{equation}
 where $\lambda_1$, $\lambda_2$, $q$, and $z$ are defined as in \eqref{ob1}. Similarly, we can   argue that \eqref{Type B} keeps hierarchy.
When $q=2$ and $z$ takes the form of \eqref{identity}, the penalty part in \eqref{Type B} become
\begin{equation}
\lambda_1\sum_{j,k}\|  \phi_{jk}\bsbx_j\odot\bsbx_k\|_2+\lambda_2\sum_j (\|b_j\bsbx_j\|_2^2+\sum_{k}\|\phi_{jk}\bsbx_j \odot \bsbx_k\|_2^2)^{\frac{1}{2}}, \label{typeB_comm}
\end{equation}
as considered by \cite{VANISH}. VANISH constructs main effects and interactions from two small sets of orthonormal basis functions in a functional regression setting. We do not pose such a restriction on the design matrix, and $p$ can be arbitrarily large.\\

The key difference between the two types of GRESH is that the penalties are imposed on the coefficients in \eqref{ob1}, but on the terms in \eqref{Type B}.
A common practice before calling a shrinkage method is normalizing/standardizing all predictors, so that it is more reasonable to use a common regularization parameter in penalizing different coefficients. In this way, \eqref{ob1} builds a model on the normalized predictors and their interactions, while Type-B amounts to forming the overall design $\breve{\bsbX}$ first and then performing the standardization.
They are not equivalent because in general,
$({\bsbx_j}/{\|\bsbx_j\|_2}) \odot  ({\bsbx_k}/{\|\bsbx_k\|_2}) \neq  ({\bsbx_j \odot \bsbx_k} )/{\|\bsbx_j \odot \bsbx_k\|_2} $.
Then, which type of GRESH is  preferable? An answer  will be given in Section \ref{sec:oracle}.

GRESH offers some general schemes for hierarchical variable selection. But it is no ordinary lasso or group lasso, since $\bsbPhi$ appears in both penalties as well as the symmetry constraint.
The main goal of this paper is to tackle some computational and theoretical challenges arising from the overlapping regularization terms in high dimensions. In computation, we would like  to develop fast and scalable algorithms (cf.\ Section \ref{sec:comp}); in theory, how to treat the penalties and the constraint \emph{jointly} to derive a sharp error bound for GRESH is intriguing and challenging (cf.\ Sections \ref{sec:oracle}, \ref{sec:minimax}).

\section{Computation}
\label{sec:comp}
 It is perhaps natural to think of using the alternating direction method of multipliers (ADMM,  cf. \cite{boyd2011distributed}) to deal with the computational challenge.   ADMM recently gains its popularity among statisticians.
In fact, \cite{bien2014convex} designed an algorithm of HL  based on      ADMM, where one of  the main ingredients is   the augmented Lagrangian
\begin{gather}
 \min_{\bsbb^\pm \in \mathbb{R}^p,\bsbPhi\in \mathbb{R}^p}\ell(\bsbb^+-\bsbb^-, \bsbPhi)+\lambda\bsb{1}^{\top}(\bsbb^+ +\bsbb^-)+\lambda\|\bsbPhi\|_1+
\langle \bsbPhi - \bsb{\Psi}, \bsb{L}\rangle +\frac{\rho}{2}\|\bsbPhi - \bsb{\Psi}\|_F^2,\notag\\  \text{s.t.} \ \bsb{\Psi}=\bsb{\Psi}^{\top},\|\bsb{\phi}_j\|_1 \leq b_j^+ + b_j^- ,b_j^+ \geq 0, b_j^- \geq 0,  \, 1\le j \le p.\label{ADMM}
\end{gather}
Here, $\bsb{L}$ is a Lagrange multiplier matrix, and $\rho>0$ is  a given constant, sometimes referred as the penalty parameter.  Although  ADMM enjoys some nice convergence properties in theory, practically  only when $\rho$ is large enough can we
obtain a  solution with good statistical accuracy. But often the larger the value of $\rho$ is, the slower the
(primal) convergence is.   For example,   the {\tt R} package  {HierNet} (version 1.6) for computing HL recommends   $\rho=n$, but for  $p=1000$, the algorithm may take several days to compute a single solution path.  There are some  empirical   schemes on how to vary $\rho$ during the iteration, but they are  ad hoc and do not always behave well.  \\

In this section, we consider a slightly more general optimization problem which includes both types of GRESH as particular instances
\begin{align}\begin{gathered}
\min_{\bsbOmega=[\bsbb,\bsbPhi^{\top}]^{\top}} \frac{1}{2}\|\bsb{y}-\bsbX\bsbb-\mdiag(  \bsbZ\bsbPhi\bsbZ^\top)\|_2^2+\|\bsblambda_{b}\odot\bsbb\|_1 \\ \qquad     +\|\bsbLambda_{\Phi}\odot\bsbPhi\|_1+ \|\bsbLambda_{\Omega}^{\top} \odot\bsbOmega^{\top}\|_{2,1}
\quad\text{s.t.}\quad \bsbPhi=\bsbPhi^{\top},
\end{gathered}\label{general0}
\end{align}
where $\bsbX, \bsbZ\in \mathbb R^{n\times p}$, and $\bsblambda_{b}$, $\bsbLambda_{\Phi}$ and $\bsbLambda_{\Omega}$ are non-negative regularization vector and matrices. Let $\bar \bsbZ=[\bsbz_1 \odot \bsbz_1,  \ldots, \bsbz_1 \odot \bsbz_p, \ldots, \bsbz_p \odot \bsbz_1, \ldots, \bsbz_p \odot \bsbz_p]$, and  $\breve \bsbZ=[\bsbx_1, \bsbz_1\odot\bsbz_1, \ldots, \bsbz_1\odot\bsbz_p, \ldots, \bsbx_p,\bsbz_p\odot\bsbz_1, \ldots, \bsbz_p\odot\bsbz_p]$.
 We   assume that   $\bsbLambda_{\Omega}=\bsb{1} \bsb{\lambda}_{\Omega}^{\top}$ for some $ \bsb{\lambda}_{\Omega}\in \mathbb R^{p}$ in developing the algorithm. (Since our algorithm applies  to a general   $\bar {\bsb{Z}}\in \mathbb R^{n\times p^2}$ that is not necessarily symmetric, this is without loss of generality.)
In \eqref{general0}, the  $\ell_1$-type penalties are imposed on  overlapping groups of variables. It is worth noting that
 the symmetry constraint  considerably complicates the grouping structure. Without it, the variable groups can be shown to follow a tree structure, for which efficient algorithms can be developed on  \cite{Jenatton11} or \cite{simon2013sparse}.

Our algorithm follows a different track than ADMM. The details are   presented   in Algorithm \ref{alg1}. Step 1 updates $\bsbXi$ and results from   a linearization-based surrogate function. Step 2 carries out    a Dykstra-like  splitting---see, e.g., \cite{Bauschke} and \cite{Sherecurrent},  by  use of two proximity   operators,    $\vec{\bsb{\Theta}}_S$ and $\Theta_S$.   Concretely, for any real number $a$,  $\Theta_S(a;\lambda)$ is given by  $\msgn(a)(|a|-\lambda)_+$ with $\msgn(\cdot)$ representing the sign function. For any vector $\bsb{a}$,   $\Theta_S(\bsb{a}; {\lambda})$ is defined componentwise and the   {multivariate} version   $\vec{\bsb{\Theta}}_S(\bsb{a};\lambda)$ is given by $\bsb{a}\Theta_S(\|\bsb{a}\|_2;\lambda)/\|\bsb{a}\|_2$ if $\bsb{a} \neq \bsb{0}$ and $\bsb{0}$ otherwise.

\begin{algorithm}[!h]
\caption{{\small The GRESH  algorithm for solving the general problem \eqref{general0} } \label{alg1}}
\begin{tabbing}
\enspace \textbf{Inputs}: \\
\enspace Data:$\bsbX$, $\bsbZ$, $\bsb{y}$.  Regularization parameters:  $\bsblambda_{b},\bsb{\Lambda}_{\Phi},\bsb{\lambda}_{\Omega}$.  \\
\enspace \textbf{Initialization}: \\
\enspace $i\gets 0$,  $\tau$ large enough (say $\tau = \|\breve \bsbZ\|_2$); \\
\enspace   $\bsblambda_{b} \leftarrow \bsblambda_{b}/\tau^2$,  $\bsbLambda_{\Phi} \leftarrow (\bsbLambda_{\Phi}+\bsbLambda_{\Phi}^{\top})/(2\tau^2)$, $\bsb{\lambda}_{\Omega} \leftarrow \bsb{\lambda}_{\Omega}/\tau^2$.\\
\textbf{ repeat} \\
\enspace 1.
 $\bsb{\Xi}_{\Phi}\gets \bsbPhi^{(i)}+\bsbZ^{\top}\mbox{diag}\{\bsb{y}-\bsbX\bsbb^{(i)}-\bar\bsbZ \vect(\bsbPhi^{(i)}))\}\bsbZ/\tau^2$,
\\
$\quad\,$ $\bsb{\Xi}_{b}\gets \bsbb^{(i)}+\bsbX^{\top}(\bsb{y}-\bsbX\bsbb^{(i)}- \bar\bsbZ \vect(\bsbPhi^{(i)}) )/\tau^2$,  $\bsb{\Xi}\leftarrow [\bsb{\Xi}_{b}, \bsb{\Xi}_{{\Phi}}^{\top}]^{\top}$; \\
\enspace 2.  $\bsb{P} \leftarrow \bsb{0}$, $\bsb{Q}\leftarrow \bsb{0}$;\\
\qquad \textbf{repeat}\\
\qquad \qquad (i) \ ${\bsbOmega}[:, k]\gets \vec{\bsb{\Theta}}_{S}(\bsb{\Xi}[:,k]+\bsb{P}[:, k];\bsb{\lambda}_{\Omega}[k]),  \forall k: 1 \leq k \leq p$;\\
\qquad \qquad (ii)\ $\bsb{P}\gets\bsb{P}+\bsb{\Xi}-\bsb{\Omega}$;\\
\qquad\qquad(iii)\ $\bsb{\Xi}[1,:]\gets\Theta_{S}(\bsb{\bsbOmega}[1,:]+\bsb{Q}[1,:];\bsblambda_{b}^{\top})$;\\
\qquad\qquad (iv) $\bsb{\Omega}[2\mbox{:end},:] \gets(\bsb{\Omega}[2\mbox{:end},:] + \bsb{\Omega}[2 \mbox{:end},:]^{\top})/2$;\\
\qquad\qquad(v) \ $\bsb{\Xi}[2\mbox{:end},:]\gets\Theta_{S}(\bsb{\Omega}[2\mbox{:}\mbox{end},:]+\bsb{Q}[2\mbox{:end},:]; \bsbLambda_{\Phi})$;\\
\qquad\qquad (vi)\ $\bsb{Q}\gets\bsb{Q}+\bsb{\Omega}-\bsb{\Xi}$;\\
\qquad\textbf{until} convergence\\
\enspace 3. $\bsb{\Omega}^{(i+1)} \leftarrow \bsb{\Xi}$;\\
\enspace 4. $i \leftarrow i+1$;\\
\enspace \textbf{until} convergence \\
\enspace \textbf{Output}
$\hat{\bsbOmega}$
\end{tabbing}
\end{algorithm}

 The GRESH algorithm is  easy to implement and involves   no complicated matrix operations  such as matrix inversion. Moreover, it does not contain  sensitive algorithmic parameters like $\rho$ in ADMM, and needs no line search.  Theorem \ref{th_conv} provides a universal theoretical choice for $\tau$ to guarantee the  global optimality of $\hat{\bsbOmega}$.
In particular,  strict {iterate} convergence,  in addition to  function-value convergence, can be established, which   is considerably stronger than an ``every accumulation point'' type conclusion in many numerical studies. For clarity, we assume  that the inner iteration runs till convergence, but this is   unnecessary;  see  {Remark 2} below.
\begin{theorem}
\label{th_conv}
Suppose $\bsblambda_{b} \geq \bsb{0}$, $\bsbLambda_{\Phi} \geq \bsb{0}$, ${\bsb{\lambda}}_{\Omega} >\bsb{0}$. For any  $\tau >\|\breve \bsbZ\|_2$ and any starting point $\bsb{\Omega}^{(0)}$, the sequence of iterates $\{\bsb{\Omega}^{(i)}\}$ converges to a {\emph{globally}} optimal solution of \eqref{general0}.
\end{theorem}

\noindent \textbf{Remark 1}. The conclusion  in the theorem holds for   ``$\ell_1+\ell_2$'' type penalties as well   \citep{elastic,Berhu}. For   the associated proximity operators, see    \cite{he2013stationary}.   In hierarchical variable selection, adding an $\ell_2$-type  shrinkage is particularly helpful to compensate for  model collinearity.

\noindent \textbf{Remark 2}. Neither the convergence of  iterates   nor the optimality guarantee  requires the full convergence of the inner loop; see  Appendix \ref{sec:appConv}
for more detail. Various stopping criteria  can be employed, e.g.,     \cite{Sch11}.   In our experience,  running  (i)--(vi) for   a few steps (say 10) usually suffices.


\noindent \textbf{Remark 3}. Algorithm \ref{alg1} and Theorem \ref{th_conv}  can be extended beyond quadratic loss functions.   When $\ell$ takes the binomial deviance in classification problems, the  first step of Algorithm 1   becomes
 $\bsb{\Xi}_{b}\gets\bsbb^{(i)}+\bsbX^{\top}(\bsb{y}-\pi(\bsbX\bsbb^{(i)}+ \bar\bsbZ\vect(\bsbPhi^{(i)})))/\tau^2$, $\bsb{\Xi}_{{\Phi}}\gets \bsb{\Phi}^{(i)}+\bsbZ^{\top}\mbox{diag}\{\bsb{y}-\pi(\bsbX\bsbb^{(i)}+\bar\bsbZ\vect(\bsbPhi^{(i)}))\}\bsbZ/\tau^2$,
 $\bsb{\Xi}\leftarrow [\bsb{\Xi}_{b}, \bsb{\Xi}_{\Phi}^{\top}]^{\top}$, where $\pi(t) = 1/(1+\exp(-t))$  and extends   componentwise to vectors. We can show that theoretically, choosing   $\tau>\|\breve \bsbZ\|_2/{2}$ guarantees the convergence of the algorithm. %

\noindent \textbf{Remark 4}. We recommend applying Nesterov's    first  acceleration  in implementations  \citep{nesterov2007gradient}. In more detail, it uses a momentum update  of $\bsb{\Xi}$ in Step 1: If $i=0$, $\bsb{\Xi}_{b}\gets \bsbb^{(i)}+\bsbX^{\top}(\bsb{y}-\bsbX\bsbb^{(i)}- \bar\bsbZ\vect(\bsbPhi^{(i)}) )/\tau^2$, $\bsb{\Xi}_{\Phi}\gets \bsbPhi^{(i)}+ \bsbZ^{\top}\mbox{diag}\{\bsb{y}-\bsbX\bsbb^{(i)}-\bar\bsbZ\vect(\bsbPhi^{(i)}))\}\bsbZ/\tau^2$; if $i>0$, $\bsb{\Xi}_{b} \gets (1-\omega_i)\bsb{\Xi}_{b}+\omega_i(\bsbb^{(i)}+\bsbX^{\top}(\bsb{y}-\bsbX\bsbb^{(i)}- \bar\bsbZ\vect(\bsbPhi^{(i)}) )/\tau^2)$, $\bsb{\Xi}_{\Phi} \gets (1-\omega_i)\bsb{\Xi}_{\Phi}+\omega_i( \bsbPhi^{(i)}+\bsbZ^{\top}\mbox{diag}\{\bsb{y}-
 \bsbX\bsbb^{(i)}-\bar\bsbZ\vect(\bsbPhi^{(i)})\}\bsbZ/\tau^2),$
where $\omega_i=(2i+3)/(i+3)$.  Empirically, the number of iterations can be reduced by about 40\%   in comparison to the non-relaxed form.

\section{Non-asymptotic analysis }
\label{sec:oracle}

In this section, given any   $\bsbDelta \in \mathbb{R}^{(p+1) \times p}$ and  $\mathcal{J}_G \subset [p]$, we use  $\bsbDelta_{\mathcal{J}_G}$ to denote the submatrix $\bsbDelta[:,\mathcal{J}_G]$.  Given any $\mathcal{J}_e \subset [p^2]$, $\|\vect(\bsbDelta_{\Phi})_{\mathcal{J}_e}\|_1$ and  $\|\vect(\bsbDelta_{\Phi})_{\mathcal{J}_e}\|_2$ are   abbreviated as $\|(\bsbDelta_{\Phi})_{\mathcal{J}_e}\|_1$  and $\|(\bsbDelta_{\Phi})_{\mathcal{J}_e}\|_2$, respectively,  when there is no ambiguity.

 In this multi-regularization setting,  some standard  treatments of the stochastic term   do not give sharp   error rates.
  In particular,     applying   $\langle \bsbeps, \allowbreak \bar{\bsbX}\vect(\bsbDelta_{\Phi}) \rangle \allowbreak \leq \|\bar{\bsbX}^{\top}\bsbeps  \|_{\infty}\|\bsbDelta_{\Phi}\|_1$  and  $\langle \bsbeps, \breve{\bsbX}\vect(\bsbDelta)\rangle \leq \|\breve{\bsbX}^{\top}\bsbeps  \|_{2,\infty}\|\bsbDelta^{\top}\|_{2,1}$   as commonly used in the    literature \citep{Bickel09,Loun2011,negahban2012,geer14}  would  yield a prediction error bound of the order  $\sigma^2(J_e \log p + J_Gp)$, which is, ironically, much worse than the error rate of  LASSO or Group LASSO (G-LASSO). Our analysis relies on   two interrelated inequalities derived from the   statistical  and computational properties of GRESH estimators. See Appendix  \ref{sec:appOracle}   for more technical detail.

First, let's    consider the  \textbf{Type-A}   problem \eqref{typeA_comm}, with $\lambda_1$, $\lambda_2$ redefined:
     \begin{equation}
  \min_{\bsbOmega=[\bsbb,\bsbPhi^{\top}]^{\top}}\frac{1}{2}\|\bsb{y}-\bsbX\bsbb-\bar{\bsbX}\vect(\bsbPhi)\|_2^2+ \lambda_1\|\breve{\bsbX}\|_2\|\bsb{\bsbPhi}\|_1  +\lambda_2\|\breve{\bsbX}\|_2\|\bsb{\bsbOmega}^{\top}\|_{2,1}  \ \text{s.t.} \ \bsbPhi=\bsbPhi^{\top}.\label{typeA_crit_oracle}
  \end{equation}
    Let   $\bsb{\hat{\Omega}}=[\hat{\bsb{b}},\hat{\bsb{\bsbPhi}}^{\top}]^{\top}$ be any global minimizer of \eqref{typeA_crit_oracle}. We are interested in its prediction accuracy measured by $M(\hat{\bsbb}-\bsbb^*,\hat{\bsbPhi}-\bsbPhi^*)$, where   \begin{equation}
  M(\bsbb,\bsbPhi)=\|\bsbX\bsbb+\bar{\bsbX}\vect(\bsbPhi)\|_2^2.
\end{equation}
The predictive learning perspective is   always legitimate  in evaluating the quality of the estimator regardless of   the signal-to-noise ratio. To guarantee small predictor errors when using   a convex method, the design matrix must satisfy certain incoherence conditions, one of the most popular     being  the restricted eigenvalue (RE) assumption   \citep{Bickel09,Loun2011}. In the following, we give an extension of RE in the hierarchy setting, with   the  restricted  cone defined with   both $\ell_1$ and group-$\ell_1$ penalties. A less intuitive  but technically much   less demanding condition is used in the proof.
   \\

\noindent {\sc Assumption $\mathcal{A}(\mathcal{J}_e,\mathcal{J}_G,\vartheta,\kappa)$}. Given $\mathcal{J}_e \subset [p^2]$,  $\mathcal{J}_G \subset [p]$, $\kappa\ge 0$ and a  constant  $\vartheta\ge 0$,  for any $\bsbDelta = [\bsbDelta_{b},\bsbDelta_{\Phi}^{\top} ]^{\top} \in \mathbb{R}^{(p+1) \times p}$ satisfying  $\bsbDelta_{\Phi}=\bsbDelta_{\Phi}^{\top}$ and  $\|(\bsbDelta_{\Phi})_{\mathcal{J}_e^c}  \|_1+\|(\bsbDelta_{\mathcal{J}_G^c})^{\top}\|_{2,1} \leq (1+\vartheta)(\|(\bsbDelta_{\Phi})_{\mathcal{J}_e}\|_1+\|(\bsbDelta_{\mathcal{J}_G})^{\top}\|_{2,1})$, the following inequality holds
  \begin{equation}
 \kappa\|\breve{\bsbX}\|_2^2(\|(\bsbDelta_{\Phi})_{\mathcal{J}_e }\|^2_2+\|\bsbDelta_{\mathcal{J}_G}\|^2_F) \leq \|\bsbX\bsbDelta_{b}+\bar{\bsbX}\vect(\bsbDelta_{\Phi})\|_2^2. \label{RE}
\end{equation}

  The rate choices of the regularization parameters play a major role in prediction. We choose  $\lambda_1$ and $\lambda_2$ in \eqref{typeA_crit_oracle} according to
\begin{equation}
  \lambda_1  = A_{1}\sigma\sqrt{\log (ep)}, \quad \lambda_2 =  A_{2}\sigma\sqrt{\log (ep)},  \label{lambdarate}
\end{equation}
where $A_1$, $A_2$ are large constants.
\eqref{lambdarate}  is  quite different from   the typical choice in      group-$\ell_1$ penalization;  see Remark 2 for more detail.

The following theorem states a  non-asymptotic oracle inequality as well as a model cardinality bound for   GRESH estimators.  For convenience,  we use  abbreviated symbols   $\hat{J}_G=J_G(\hat \bsbb, \hat \bsbPhi)$ and  $\hat{J}_e=J_e(\hat\bsbPhi)$ for the estimate, and ${J}_G=J_G( \bsbb,  \bsbPhi)$ and  ${J}_e=J_e(\bsbPhi)$ for the reference signal.

\begin{theorem}
  \label{th_oracle}
  Assume $\boldsymbol{\epsilon}\sim N(\bsb{0}, \sigma^2 \boldsymbol{I})$.   Let $\hat{\bsbOmega}=\big[\hat{\bsb{b}},\hat{\bsbPhi}^{\top}\big]^{\top}$ be a global minimizer of \eqref{typeA_crit_oracle}. Then under \eqref{lambdarate},   for any sufficiently large constants $A_1, A_2$, the following oracle inequality holds for any $(\bsbb , \bsbPhi)\in \mathbb R^{p}\times \mathbb R^{p\times p}$
  \begin{equation}
  \begin{aligned}
  \mathbb{E}[M(\hat{\bsbb}-\bsbb^*,\hat{\bsbPhi}-\bsbPhi^*)] \lesssim M(\bsbb-\bsbb^*,\bsbPhi-\bsbPhi^*)+(1\vee \frac{1}{\kappa})\sigma^2 (J_e+J_G)\log p+ \sigma^2,
   \label{predoracle}
  \end{aligned}
  \end{equation}
 provided that $(\bsbX, \bar{\bsbX},\breve{\bsbX})$ satisfies $\mathcal{A}(\mathcal{J}_e,\mathcal{J}_G,\vartheta,\kappa)$ for some    $\kappa>0$ and some constant $\vartheta\ge 0$. Furthermore, under the same regularity condition, the overall sparsity of the obtained model is  controlled by
  \begin{equation}
  \begin{aligned}
  \mathbb{E}[\hat{J}_e] +  \mathbb{E}[\hat{J}_G] \lesssim \{ M(\bsbb^*-\bsbb,\bsbPhi^*-\bsbPhi)+\sigma^2\}/\{\sigma^2\log(ep)\}+ J_e+ J_G. \label{cardbound}
  \end{aligned}
  \end{equation}
\end{theorem}

\noindent \textbf{Remark 1}.   Letting $\bsbb = \bsbb^*$ and $\bsbPhi=\bsbPhi^*$ in \eqref{predoracle}, we obtain an error bound no larger than   $\sigma^2(J^*_e+J^*_G)\log p$ (omitting constant factors). This indicates that GRESH not only guarantees SH, but can give an error rate as low as that of LASSO. The existence of the bias term  $M(\bsbb-\bsbb^*, \bsbPhi-\bsbPhi^*)$  makes our results applicable to {approximately} sparse signals,  which is of   practical significance.
The theorem does not require the spectral norms of the design matrices $\bsbX$, $\bar{\bsbX}$ and $\breve{\bsbX}$ to be bounded above by  $O(\sqrt{n})$ as   assumed in, for example,   \cite{zhang2008sparsity} and \cite{Bickel09}. In addition,   the true signal $\bsbOmega^*$ and the reference signal $\bsbOmega$   in the theorem need not    obey SH.

\noindent\textbf{Remark 2}.   It is widely acknowledged that the penalty parameter for a grouped $\ell_1$ penalty should be adjusted by the group size \citep{grouplasso}.  In fact,    $\lambda_2$ would be of order  $\sigma\sqrt{p+\log p}$ from \cite{Loun2011} and \cite{wei2010consistent},   in light of the fact that there are  $p$ groups of size $(p+1)$ in $\|\bsbOmega\|_{2,1}$.
 Perhaps surprisingly, this parameter choice  becomes suboptimal in hierarchical variable selection. In fact, due to the presence of multiple penalties, we show in the proof that   \eqref{lambdarate}  suffices to suppress   the noise, which in turn  leads to a reduced  error rate. Such a novel finding is  owing to      the careful treatment of the stochastic term,    
 which     is generally applicable to    overlap group lasso \citep{Jenatton11}. The conclusion that $\lambda_1$ and $\lambda_2$ are of essentially the same rate also facilitates parameter tuning, since  one just needs to search along a one-dimensional grid.

\noindent\textbf{Remark 3}. Theorem \ref{th_oracle} can be extended to sub-Gaussian  $\vect(\bsbeps)$  with mean 0 and  its $\psi_2$-norm bounded by $\sigma$, which covers more     noise distributions.  High probability form results of the prediction error  can be obtained as well:    \eqref{predoracle} and \eqref{cardbound}, without the expectation and  the additive $\sigma^2$ term, hold with probability at least $1-C p^{-c \min\{A_1^2,   A_2^2\}}$ for some universal constants $C$ and $c$, and so       $\hat{J}_e+\hat{J}_G \lesssim J_e^* + J_G^*$  with high probability. Moreover, in Appendix \ref{sec:coord},
we   show how to adapt our proof to deliver a coordinatewise error bound which can be used for recovering the sparsity pattern of the true signal.

\noindent\textbf{Remark 4}. For the WH version of \eqref{typeA_crit_oracle} (without the symmetry condition), following the lines of the proof of Theorem \ref{th_oracle}, we can show its error rate  is of the order    $\sigma^2 \{J_e^{w}(\bsbPhi)+J_G^{w}(\bsbOmega)\}\log p$, where    $J_G^{w}(\bsbOmega)=J_G(\bsbOmega')$,   $J_e^{w}(\bsbPhi)=J_e(\bsbPhi')$, with $\bsbOmega' = [\bsbb, \bsbPhi'^\top]^\top$ and $\phi_{k j}'=\phi_{k j}+\phi_{j k}$ for $k\geq j$ and $0$ otherwise. The associated regularity condition uses  $\mathcal J_e^w$, $\mathcal J_G^w$, in place of   $\mathcal J_e$, $\mathcal J_G$, respectively,   and  does not require $\bsbDelta_{\Phi}$ to be symmetric. Details are not reported in the paper.
\\

Similarly, we can derive an oracle inequality for GRESH estimators of Type-B.   Let $\bsbX^s$ be the column-scaled $\bsbX$  such that  the $\ell_2$-norm of each of its columns equals 1.   $\bar{\bsbX}^s $ and $\breve{\bsbX}^s $ are similarly defined.  The corresponding coefficients,  denoted by $\bsbb^s$ and $\bsbPhi^s$, satisfy $\bsbX\bsbb =\bsbX^s\bsbb^s$, $\bar{\bsbX}\vect(\bsbPhi)=\bar{\bsbX}^s\vect(\bsbPhi^s)$. Let $(\hat{\bsb{b}}^s, \hat{\bsb{\bsbPhi}}^{s})$ be a global minimizer of the  {scaled} \textbf{Type-B} problem:
$\min_{\bsbOmega^s=[\bsbb^s,\bsbPhi^{s\top}]^{\top}}\frac{1}{2}\|\bsb{y}-\bsbX^s\bsbb^s-\bar{\bsbX}^s\vect(\bsbPhi^s)\|_2^2 + \lambda_1\|\breve{\bsbX}^s\|_2\|\bsb{\bsbPhi}^s\|_1 +\lambda_2\|\breve{\bsbX}^s\|_2\|(\bsbOmega^s)^{\top}\|_{2,1} $ s.t. $\bsbPhi^s=(\bsbPhi^s)^{\top}$. As aforementioned, the problem can not be reduced to  \eqref{typeA_crit_oracle}  because $\bar{\bsbX}^s[:,jk]$ does not equal    $\bsbX^s[:,j]\odot\bsbX^s[:,k]$ in general.  \\

 \noindent {\sc Assumption $\mathcal{A'}(\mathcal{J}_e,\mathcal{J}_G,\vartheta,\kappa')$}. Given $\mathcal{J}_e \subset [p^2]$,  $\mathcal{J}_G \subset [p]$, and positive constants $\kappa'$ and $\vartheta$,   for any $\bsbDelta = [\bsbDelta_{b},\bsbDelta_{\Phi}^{\top} ]^{\top}$ satisfying $\bsbDelta_{\Phi}=\bsbDelta_{\Phi}^{\top}$ and  $\|(\bsbDelta_{\Phi})_{\mathcal{J}_e^c}  \|_1+\|(\bsbDelta_{\mathcal{J}_G^c})^{\top}\|_{2,1} \leq (1+\vartheta)(\|(\bsbDelta_{\Phi})_{\mathcal{J}_e}\|_1+\|(\bsbDelta_{\mathcal{J}_G})^{\top}\|_{2,1})$, the following inequality holds
  \[
 \kappa'\|\breve{\bsbX}^s\|_2^2(\|(\bsbDelta_{\Phi})_{\mathcal{J}_e}\|^2_2+\|\bsbDelta_{\mathcal{J}_G}\|^2_F) \leq \|\bsbX^s\bsbDelta_{b}+\bar{\bsbX}^s\vect(\bsbDelta_{\Phi})\|_2^2.
 \]

\noindent \textbf{Theorem 2'} \emph{Under the same conditions as in Theorem \ref{th_oracle} and  with  $\mathcal{A'}(\mathcal{J}_e,\mathcal{J}_G,\vartheta,\delta'_{\mathcal{J}_e,\mathcal{J}_G})$ in place of $\mathcal{A}(\mathcal{J}_e,\mathcal{J}_G,\vartheta,\kappa)$, \eqref{predoracle} and \eqref{cardbound} hold}. \\

The error bounds of the  two types of GRESH are of  the same order, but their regularity conditions place different requirements on the design. We   performed extensive simulation studies to compare  $\mathcal{A}$ and $ \mathcal{A'}$, and found that for the same $\vartheta$, $\kappa< \kappa'$  usually holds, which suggests    the penalization on the basis of  terms seems more appropriate  than that on the  coefficients. Therefore, we recommend Type-B regularization for hierarchical variable selection.

\section{Minimax lower bound and error rate comparison}
\label{sec:minimax}
  In this section, we show that in a minimax sense, the error rate we obtained in Theorem \ref{th_oracle} is  minimax optimal up to some logarithmic factors.
Consider two signal classes having  hierarchy and joint sparsity:
 \begin{gather}
\text{\textbf{SH}}(J_G, J_e)= \{\bsbOmega=[\bsbb,\bsbPhi^{\top}]^{\top}: \bsbOmega  \text{ obeys {SH}}, \bsbPhi=\bsbPhi^\top, J_G(\bsbOmega) \leq J_G, J_e(\bsbPhi)\leq J_e   \},\\
\text{\textbf{WH}}(J_G, J_e)= \{\bsbOmega=[\bsbb,\bsbPhi^{\top}]^{\top}: \bsbOmega  \text{ obeys {WH}},  J_G^{w}(\bsbOmega) \leq J_G, J_e^w(\bsbPhi)\leq J_e   \},
\end{gather}
where ${1} \leq J_G \leq p$, ${1} \leq J_e \leq p J_G$.
 Recall the definitions of  $J_G^{w}$ and $J_e^w$ in  {Remark 4} following Theorem \ref{th_oracle}.
Let $\ell(\cdot)$ be a nondecreasing loss function with $\ell(0)=0$, $\ell \not\equiv0$. Under   some regularity  assumptions, we study the minimax lower bounds for strong and weak hierarchical variable selection.
\\

\noindent {\sc Assumption} $\mathcal{B^S}(J_G, J_e)$.
     For any $\bsbOmega=[\bsbb,\bsbPhi^{\top}]^{\top} \in \mathbb{R}^{(p+1) \times p}$ satisfying that  $\bsbPhi$  is symmetric,  $J_e(\bsbPhi)\leq J_e$ and $J_G(\bsbOmega) \leq J_G $, $\underline{\kappa} \|\bsbOmega\|_F^2 \leq  \|\breve{\bsbX}\vect(\bsbOmega)\|_2^2 \leq \overline{\kappa} \|\bsbOmega\|_F^2
    $ holds, where    $\underline{\kappa}/\overline{\kappa}$ is a positive constant.
\\

\noindent {\sc Assumption} $\mathcal{B^W}(J_G, J_e)$.
     For any $\bsbOmega=[\bsbb,\bsbPhi^{\top}]^{\top} \in \mathbb{R}^{(p+1) \times p}$ satisfying that  $J_e^{w}(\bsbPhi)\leq J_e$ and $J_G^{w}(\bsbOmega) \leq J_G $, $\underline{\kappa} \|\bsbOmega\|_F^2 \leq  \|\breve{\bsbX}\vect(\bsbOmega)\|_2^2 \leq \overline{\kappa} \|\bsbOmega\|_F^2
    $ holds, where     $\underline{\kappa}/\overline{\kappa}$ is a positive constant.

\begin{theorem}
\label{th_minimax}
(i) Strong hierarchy. Assume $\bsb{y}=\bsbX\bsbb^* + \bar{\bsbX}\vect(\bsbPhi^*) +\bsbeps$ with $\bsbeps \sim N(\bsb{0},\sigma^2\bsb{I})$, $J_G\geq 2$, $p\geq 2$, $J_e\geq 1$, $n\geq 1$, $J_G\le p/2$, $J_e\leq J_G^2 /2$,
   and $\mathcal{B^{S}}(2J_G, 2J_e)$ is satisfied.  Then there exist positive constants $C$, $c$ (depending on $\ell(\cdot)$ only) such that
\begin{equation}
    \inf_{\hat{\bsbOmega}}\sup_{\bsbOmega^* \in \text{SH}(J_G,J_e)}\mathbb{E}[\ell(M(\hat{\bsbb}-\bsbb^*,\hat{\bsbPhi}-\bsbPhi^*)/(C P_o(J_e,J_G)))] \geq c>0, \label{minimaxlowerbound}
\end{equation}
where   $\bsb{\hat{\Omega}}$ denotes any  estimator, and
\begin{align}
P_o(J_e,J_G)= \sigma^2\{J_e\log(eJ_G^{2}/J_e) + J_G\log(ep/J_G)\}.
\end{align}
(ii) Weak hierarchy. Let $J_G \geq 1$, $J_e \geq 1$, $n\geq 1$, $p\geq 2$, $J_G\le p/2$, $J_e\leq J_G p /2$.
Under the same model assumption and $\mathcal{{B}^{W}}(2J_G, 2J_e)$,
\eqref{minimaxlowerbound}  holds if   $\text{SH}(J_G,J_e)$ is replaced by $ \text{WH}(J_G,J_e)$ and  $P_o$ is replaced by
\begin{align}
P_o'(J_e,J_G)= \sigma^2\{J_e\log(eJ_Gp/J_e) + J_G\log(ep/J_G)\}.
\end{align}
\end{theorem}

  We give some examples of $\ell$ to illustrate the conclusion. For SH, using the indicator function $\ell(u)=1_{u\geq 1}$, we know that for any estimator $(\hat{\bsbb},\hat{\bsbPhi})$,
\begin{equation*}
M(\hat{\bsbb}-\bsbb^*,\hat{\bsbPhi}-\bsbPhi^*)  \gtrsim \sigma^2\left(J_e\log(eJ_G^2/J_e) + J_G\log(ep/J_G)\right)
\end{equation*}
occurs with positive probability,  under some mild conditions. For $ \ell(u)=u$, Theorem \ref{th_minimax} shows that the risk $\mathbb{E}[M(\hat{\bsbb}-\bsbb^*,\hat{\bsbPhi}-\bsbPhi^*)]$ is bounded from below by $P_o(J_e,J_G)$   up to some multiplicative constant. Because $J_G\leq J_e \leq J_G p$ and  $J_G \leq p$, it is easy to see that the minimax rates are no larger than  the error rate obtained  in Theorem \ref{th_oracle}.
\\

A comparison of some popular methods follows, where we can   see the  benefits of  hierarchical variable selection.  In our context, LASSO  solves
 $\min_{\bsbb, \bsbPhi:\bsbPhi=\bsbPhi^{\top}}\allowbreak
  \ell(\bsbb, \bsbPhi)+\lambda(\|\bsbb\|_1 +\|\bsbPhi\|_1)$.    From \cite{Bickel09},  the estimator has a prediction error  of the order $\sigma^2 (J_e^{}+J^{10}+J^{11})\log p$.  G-LASSO, with the optimization problem defined by $\min_{\bsbb, \bsbPhi:\bsbPhi=\bsbPhi^{\top}}\ell(\bsbb, \bsbPhi)+\lambda\|\bsbOmega^{\top}\|_{2,1}$,      automatically maintains SH and has an   error rate of    $\sigma^2J_G p$ \citep{Loun2011}.
In general, there is no clear winner between the two. Let's turn to a particularly interesting    {case where }   $J^*_e \ll J^*_G p$, i.e., the existence of a main effect in the model does not indicate {that all} its associated interactions must be relevant. In this scenario, LASSO always outperforms G-LASSO, although it   does not possess the SH property.   GRESH achieves the same low error rate and guarantees hierarchy, because under SH, $J^{11}+J^{10}=J_G$.

The error rate proved  in \eqref{predoracle}  does not always beat that of G-LASSO,   because only   large values of $A$  are considered in Theorem \ref{th_oracle}. Yet, even in the worst case when $J^{*}_e \asymp J^*_Gp $,  GRESH   is only a logarithmic factor worse. In practical data analysis,  there will  {be no} performance loss, because when $\lambda_1=0$, GRESH degenerates to G-LASSO.

\begin{table}[h]
 \centering  \small{
\caption{\label{tab:ratecomp} {Error rate comparison between LASSO, G-LASSO, and GRESH, where     $\sigma^2$ and other constant factors are omitted.}}%
\vspace{.15in}
 \begin{tabular}{lc}
 \hline
 LASSO    &$(J_{e}+J^{11}+J^{10})\log p$         \\       
 G-LASSO   &$J_G p$\\
 GRESH      &$(J_G+J_e)\log p$\\
 \hline
 \textbf{Minimax}:  &$J_G\log(ep/J_{G})+J_e\log (e J_G^2 /J_e)$ \\
 \hline
 \end{tabular}
} \end{table}

\section{Experiments}
\label{sec:exp}
\subsection{Simulations}
\label{subsec:exp_simu}
   In this part, we perform some simulation studies to compare the performance of HL and GRESH (of Type B, cf. \eqref{typeB_comm})  in terms of prediction accuracy,  selection consistency, and computational efficiency. We use a Toeplitz design to generate all main predictors, with the correlation between $\bsb{x}_i$ and $\bsb{x}_j$  given by $0.5^{|i-j|}$. The true coefficients  $\bsbb^*$ and $\bsbPhi^*$ (symmetric) are generated according to the following three setups.\\

\noindent \textbf{Example 1.}  \ $n=40$, $p=100$ or  $200$ (and so $p + {p \choose 2} =5050$ or  $20100$).  $\bsbb^*=[3,1.5,0,0,2, 2, 0,\cdots,0]^{\top}, \bsbPhi^*=\bsb{0}$, $\sigma^2=1$. No interactions are relevant to the response variable.  SH is satisfied.

 \noindent \textbf{Example 2.}  \ $n=150$, $p=50$ or $100$ (and so  $p + {p \choose 2}=1275$ or $5050$).   $\bsbb^*=[3,3,3,3,3,3,3,3,3,0, \cdots,0]^{\top}$, $\bsbPhi^*[1,2]=\bsbPhi^*[1,3]=\bsbPhi^*[4,5]=\bsbPhi^*[4,6]=\bsbPhi^*[7,8]=\bsbPhi^*[7,9]=3$, $\sigma^2=1$. The model involves both main and interaction effects and obeys SH.

  \noindent \textbf{Example 3.}  \ $n=100$, $p=50$ or $100$ (and so $p + {p \choose 2}=1275$ or $5050$).  $\bsbb^*=[1, 1, 1, 1, 1, 1,  1, 0, \cdots, 0]^\top$, $\bsbPhi^*[i, j]=5$, $1\leq i, j\leq 5$, $i\neq j$, $\bsbPhi^*[4,5]=\bsbPhi^*[4,6]=\bsbPhi^*[4,7]=5$, $\sigma^2=1$. The true model does not have very strong main effects but satisfies  SH. \\

 All regularization parameters are tuned on a (separate) large validation dataset containing 10K observations. There is no need to perform  a full two-dimensional  grid search to find the optimal   parameters in GRESH. Rather,  motivated by Theorem \ref{th_oracle}, we set $\lambda_2=c\lambda_1$, and chose   $c=0.5$ according to    experience.    Because of the convex nature of the problem, pathwise  computation with warm starts is used.
After   variable selection, a   ridge regression model is always refitted to be used for  prediction. The official  {\tt R} package for HL is {HierNet}, implemented in {\tt C}.  We set  {\tt strong=TRUE} and    post-calibrate  HL by a restricted ridge refitting, which substantially enhances its accuracy. To make a fair comparison between HL and GRESH, we use the same error tolerance (1e-5) and the same number of grid values ($20$).
All other  algorithmic parameters in HierNet are set to their default values.
Given each setup, we repeat the experiment for 50 times and evaluate the performance of each algorithm according to the measures defined below. The test error (Err) is  the mean squared error between the true mean of $\bsb{y}$ and its estimate;  
 for robustness and stability, we report  the median test error from all runs. The joint detection (JD) rate is the fraction of $|\lbrace (i,j) : \Omega_{ij}^* \neq 0 \rbrace | \subseteq  |\lbrace (i,j): \hat{\Omega}_{ij}\neq 0\rbrace|$ among all experiments. The  missing (M) rate and the false alarm (FA) rate are the mean of $|\lbrace (i,j): \Omega_{ij}^* \neq 0, \hat{\Omega}_{ij} = 0 \rbrace| / | \lbrace (i,j) : \Omega_{ij}^* \neq 0 \rbrace|$ and the mean  $|\lbrace (i,j):\Omega_{ij}^{*}=0,  \hat{\Omega}_{ij}\neq 0 \rbrace|/|\lbrace(i,j):\Omega_{ij}^*=0\rbrace|$, respectively.  The path computational cost is the average running time of an algorithm in seconds. All the experiments were run on a PC with 3.2GHz CPU, 32GB memory and 64-bit Windows 8.1. Table \ref{tab:statperf} and Table \ref{tab:compperf}  summarize the statistical and computational results.

\begin{table}[h]
\centering
\small{
\caption{Statistical performance of HL and GRESH, measured in test error, joint detection  rate, missing rate, and false alarm rate on simulation data. All numbers are multiplied by 100. \label{tab:statperf} }
\vspace{.15in}
\setlength{\tabcolsep}{1.8mm}
 \begin{tabular}{l cccc cccc cccc}
         \hline
          \multirow{2}{*}{} & \multicolumn{4}{c}{\textbf{Ex 1} ($p=100$)} & \multicolumn{4}{c}{\textbf{Ex 2} ($p=50$)}  & \multicolumn{4}{c}{\textbf{Ex 3} ($p=50$)} \\  \cmidrule(r){2-5} \cmidrule(r){6-9} \cmidrule(r){10-13}
         & Err& JD & M & FA  &Err & JD & M & FA & Err & JD &M &FA \\ \hline

{HL} & 13.7  & 100& 0.00 & 0.00 & 14.2 & 95  & 0.24 & 0.31  & 68.0 & 90 & 0.65 & 2.71      \\
{GRESH}  &13.6 & 100  & 0.00 & 0.00  &11.6 & 100 & 0.00 &0.12  & 26.2 & 100 & 0.00 & 0.23   \\

\hline

\hline
          \multirow{2}{*}{} & \multicolumn{4}{c}{\textbf{Ex 1} ($p=200$)} & \multicolumn{4}{c}{\textbf{Ex 2} ($p=100$)}  &\multicolumn{4}{c}{\textbf{Ex 3}  ($p=100$)} \\
         \cmidrule(r){2-5} \cmidrule(r){6-9}\cmidrule(r){10-13}
         & Err & JD & M  & FA  &Err &JD &M & FA &Err & JD & M & FA \\ \hline
          {HL} &13.8 & 100 & 0.00 & 0.00  & 17.3 &90  & 0.48 &0.11 & 92.2 & 15 &24.03 &1.44   \\
{GRESH} &13.5 & 100 & 0.00 & 0.00 & 14.3 & 100 & 0.00 & 0.05 &26.9 & 100 & 0.00 & 0.08  \\

\hline
\end{tabular}
}
\end{table}

\begin{table}[h]
\centering
\small{
\caption{Path computation costs of GRESH and HL when $p = 200$ and $1000$. The computational times are in seconds unless otherwise specified.\label{tab:compperf}}
\vspace{.15in}
\setlength{\tabcolsep}{2mm}
 \begin{tabular}{l ccc cccc }
         \hline

          \multirow{2}{*}{} & \multicolumn{3}{c}{$p=200$} & \multicolumn{3}{c}{$p=1000$}  \\  \cmidrule(r){2-4} \cmidrule(l){5-7} 
        & Ex1 & Ex2  & Ex3 & Ex1  & Ex2 & Ex3  \\    \cmidrule(r){1-4} \cmidrule(l){5-7}

{HL}  &  574& 3057 & 2.9\,{hours} &  7.6\,{hours}    & \textbf{---} & \textbf{---}\\
 {GRESH}& 110 &158&128  & 1066 
          & 2.5\,hours   &3194  \\

\hline
\end{tabular}
}
\end{table}


  From Table \ref{tab:statperf},   GRESH and HL behaved equally well in  Example 1, the model of which   contains main effects only, but  GRESH is   faster. In Example 2  and Example 3, the two methods show  more differences;  see their test errors   and joint identification rates, for example.    We also noticed  that    GRESH   often gave  a more parsimonious model.  When   the  main effects are weak as in   Example 3,     HL  may miss some   genuine interaction effects.   Overall, GRESH showed comparable or better test errors. In fact, this is observed even when SH is not satisfied (results not shown in the table). We suspect that the performance differences between HL and GRESH largely result from the fact  that HL   compares $|b_j|$ with  $\|\bsb{\phi}_j\|_1$, the $\ell_1$-norm of the  {overall}   $\bsb{\phi}_j$, to realize SH, while   \eqref{typeB_comm}  groups   $b_j\bsbx_j$, $\phi_{j 1}\bsbx_j \odot \bsbx_1$, $\ldots$, and  $\phi_{j p}\bsbx_j \odot \bsbx_p$, on the term basis, to select main effects.

The {computational times} in Table \ref{tab:compperf} show the scalability of each algorithm as  $p$ varies.  When $p=1000$,  there are $500500$ variables in total, and so  HL   became    computationally prohibitive,  also   evidenced by  \cite{lim2013learning}. GRESH offered impressive computation gains  in the experiment.

\subsection{Comparison with ADMM}

This part shows that directly applying ADMM does not give a scalable algorithm for solving the optimization problem \eqref{general0} which has a large number of groups  with large group size. The detailed algorithm design is given  in  Appendix \ref{appadmm}.  We set   $\rho = 1$ in ADMM and compared it to Algorithm \ref{alg1}. The results are reported in Table \ref{tab:statperf_ADMM} and Table \ref{tab:compperf_ADMM}.  The    statistical performances of the two algorithms are close.   This is reasonable because they   solve the same optimization problem.  However,  ADMM is much slower. In the experiments,  ADMM became   infeasible  when $p=200$ or larger. 

\begin{table}[h]
\centering
\small{
\caption{Statistical performance of GRESH and ADMM, measured in test error, joint detection  rate, missing rate, and false alarm rate.  \label{tab:statperf_ADMM} }
\vspace{.15in}
\setlength{\tabcolsep}{1.8mm}

 \begin{tabular}{l cccc cccc cccc}
         \hline
          \multirow{2}{*}{} & \multicolumn{4}{c}{\textbf{Ex 1} ($p=50$)} & \multicolumn{4}{c}{\textbf{Ex 2} ($p=50$)}  & \multicolumn{4}{c}{\textbf{Ex 3} ($p=50$)} \\  \cmidrule(r){2-5} \cmidrule(r){6-9} \cmidrule(r){10-13}
         & Err& JD & M & FA  &Err & JD & M & FA & Err & JD &M &FA \\ \hline

{GRESH}  & 11.6 & 100 & 0.00 & 0.02    &11.6 & 100 & 0.00 &0.12  & 26.2 & 100 & 0.00 & 0.23   \\
ADMM & 11.7 & 100 & 0.00 & 0.02  &11.9 & 100 & 0.00 &0.12  & 26.3 & 100 & 0.00 & 0.22\\

\hline

\hline
          \multirow{2}{*}{} & \multicolumn{4}{c}{\textbf{Ex 1} ($p=100$)} & \multicolumn{4}{c}{\textbf{Ex 2} ($p=100$)}  &\multicolumn{4}{c}{\textbf{Ex 3}  ($p=100$)} \\
         \cmidrule(r){2-5} \cmidrule(r){6-9}\cmidrule(r){10-13}
         & Err & JD & M  & FA  &Err &JD &M & FA &Err & JD & M & FA \\ \hline
 {GRESH} &13.6 & 100  & 0.00 & 0.00  & 14.3 & 100 & 0.00 & 0.05 &26.9 & 100 & 0.00 & 0.08  \\
          ADMM & 16.8  & 100  & 0.00 & 0.01 & 14.5 & 100  & 0.00  &  0.06  & 27.1  & 100 & 0.00 & 0.09 \\

\hline
\end{tabular}
}
\end{table}

\begin{table}[h]
\centering
\small{
\caption{Path computation costs of GRESH and ADMM.  \label{tab:compperf_ADMM}}
\vspace{.15in}
\setlength{\tabcolsep}{1.8mm}
 \begin{tabular}{l ccc cccc  }
         \hline

          \multirow{2}{*}{} & \multicolumn{3}{c}{$p=50$} & \multicolumn{3}{c}{$p=100$}  \\  \cmidrule(r){2-4} \cmidrule(l){5-7} 
         & Ex1 & Ex2  & Ex3 & Ex1  & Ex2 & Ex3  \\    \cmidrule(r){1-4} \cmidrule(l){5-7}


{GRESH} & 8 & 11 & 9  & 27 
          & 36   &30  \\
ADMM  & 66 & 70  & 60  & 1027  & 1571  &  1439 \\

\hline
\end{tabular}
}
\end{table}

\subsection{Real data example}
\label{subsec:exp_real}
We performed hierarchical variable selection on the California housing data \citep{KelleyPace1997}. 
 The dataset consists of 9 summary characteristics for 20640 neighborhoods in California.  The response variable is the  median house value in each neighborhood. Following  \cite{ESL2}, we obtained eight household-related predictor variables:   median income, housing median age, average number of rooms and bedrooms per household,  population, average occupancy (population/households), latitude, and longitude,  denoted by {\tt MedInc},   {\tt Age}, {\tt AvgRms}, {\tt AvgBdrms},   {\tt Popu},  {\tt AvgOccu}, {\tt Lat}, and {\tt Long}, respectively.
Similar to
  \cite{Spam} and  \cite{VANISH},  $50$ nuisance features  generated as standard Gaussian random variables  were added, to make the  problem more challenging. The full quadratic model on this enlarged  dataset contains 3422 unknowns.

To prevent from getting over-optimistic error   estimates, we used a {hierarchical} cross-validation procedure where an outer 10-fold cross-validation (CV) is   for performance evaluation and the inner 10-fold selective CVs \citep{Group_She} are for parameter tuning.
We managed to run both HL and GRESH for hierarchical variable selection, with the estimates  post-calibrated by a local ridge fitting as described in Section \ref{subsec:exp_simu}. It  took us approximately one and half days to complete the  CV experiment for HL, and   about $1.6$ hours for   GRESH.
  The median and mean  test errors of the models obtained by HL  are  $530.8$ and $553.5$, respectively, and  the average  number of  selected variables is   $31.2$.   GRESH  gave   $516.9$ and $521.1$ for the median and mean test errors, ,  respectively, and selected      $17.1$ variables on average, about half of the model size of HL.

\begin{figure}[h!]
\begin{center}
$\begin{array}{cc}
\includegraphics[width=2.5in]{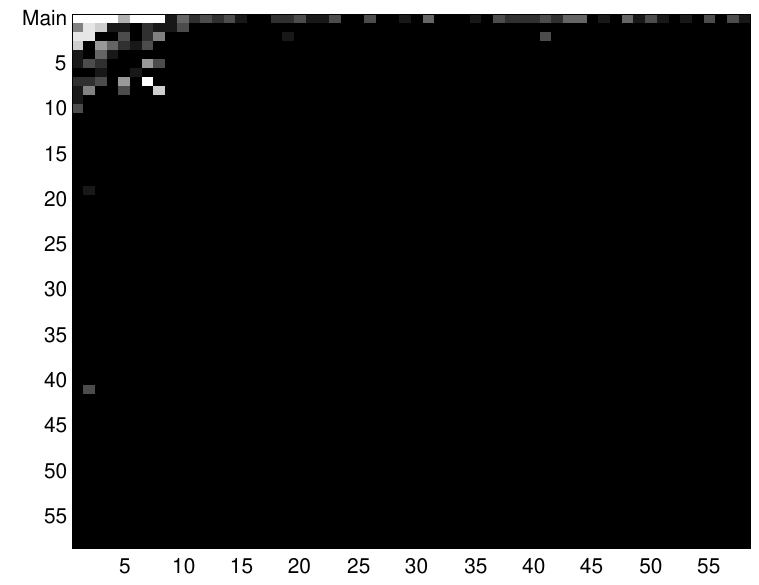}
&\includegraphics[width=2.5in]{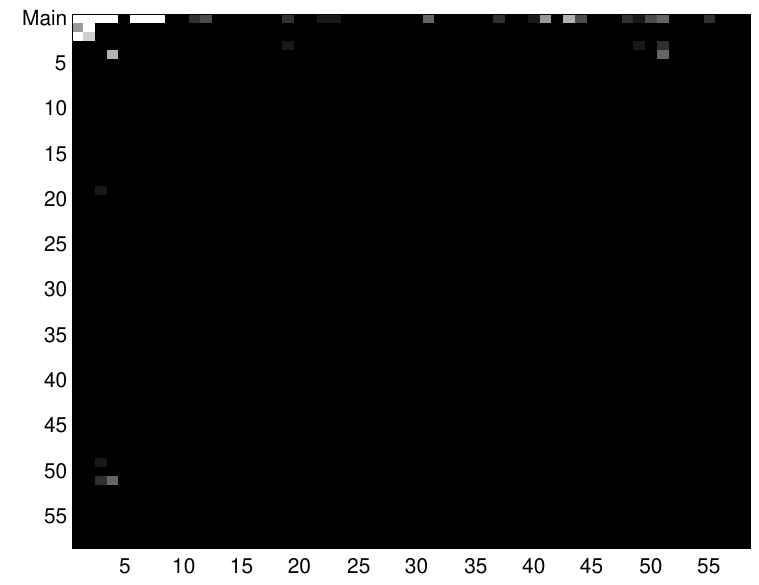}\\
\includegraphics[width=2.375in]{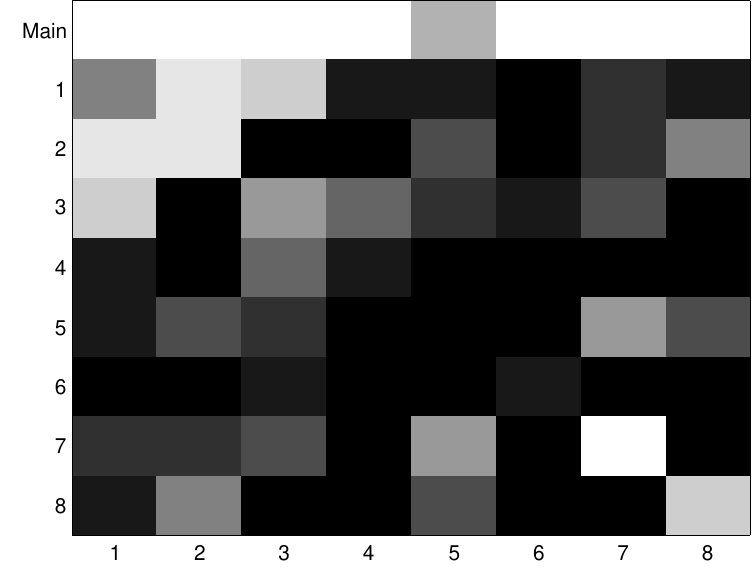}
&\includegraphics[width=2.375in]{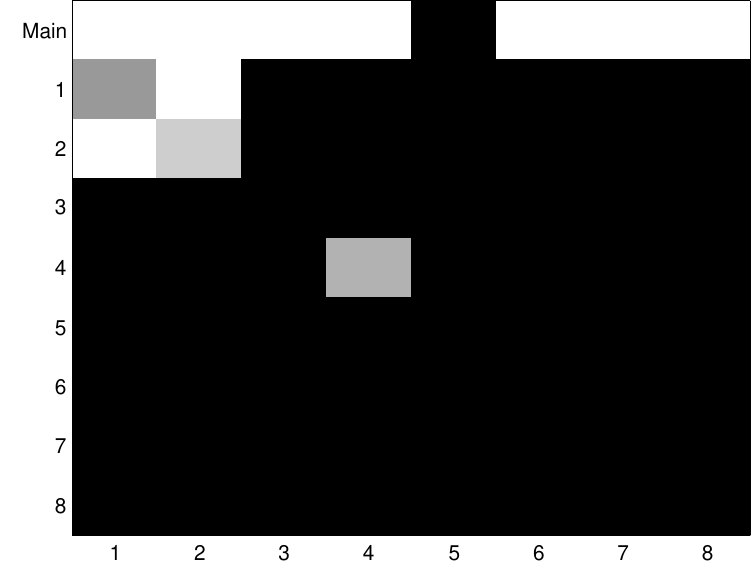}
\end{array}
$\end{center}
\caption{\small{Top panel: heat maps of  HL (left) and GRESH (right) on California housing data. Bottom panel:   heat maps of HL (left) and GRESH (right) restricted to the original 8 variables and their interactions. }}
\label{figheatmaps}
\end{figure}

To help the reader get an intuition of the selection frequencies of all predictors, we display heat maps in  Figure \ref{figheatmaps}. The two heat maps in the top panel include    all variables, and the bottom panel only shows   the heat maps restricted to the   original covariates and their interactions.
According to the figure, both methods successfully removed most of the artificially added noisy features. On average, only  $9.3$  nuisance covariates exist in the models obtained by HL, and  $5.4$  in the GRESH models. The heat maps of GRESH are  however neater. The HL selection results are less parsimonious and are perhaps more  difficult to interpret. The nonlinear terms in GRESH   include  the interaction between  {\tt MedInc} and {\tt Age},  in addition to  the quadratic effects  of {\tt MedInc}, {\tt Age},  and {\tt AvgBdrms}.         {\tt Popu} and all its associated interaction terms never got selected by GRESH. The insignificance   of {\tt Popu} can be confirmed by     more  e·lab·o·rate analysis based on gradient boosting     \citep{ESL2}.  


\appendix
\appendixpage
\section{Proofs of main theorems}
\subsection{Proof of Theorem \ref{th_conv}}
\label{sec:appConv}
Denote  the object function of \eqref{general0} by  $f(\bsb{\Omega};\bsblambda_{b},\bsbLambda_{\Phi},\bsblambda_{\Omega})=\ell(\bsbOmega)+P_1(\bsb{\Omega}; \bsblambda_{b},\\ \bsbLambda_{\Phi} )+P_2(\bsb{\Omega};\bsblambda_{\Omega})$,
where $P_1(\bsb{\Omega}; \bsblambda_{b},\bsbLambda_{\Phi} )=\|\bsblambda_{b} \odot \bsbb\|_1 + \|\bsbLambda_{\Phi} \odot \bsbPhi\|_1$, and $P_2(\bsb{\Omega};\bsblambda_{\Omega})=\|(\bsblambda_{\Omega} \bsb{1}^\top) \odot \bsbOmega^{\top}\|_{2,1}$.   Define $g(\bsbOmega,\bsbOmega'; \bsblambda_{b},\bsbLambda_{\Phi},\bsbLambda_{\Omega})
=\ell(\bsbOmega')+\frac{\tau^2}{2}\|\bsbOmega'-\bsbOmega \|_F^2+\langle \bsbX^{\top}(\bsb{y}-\bsbX\bsbb'- \bar\bsbZ \vect(\bsbPhi') ), \bsbb- \bsbb' \rangle +\langle\bsbZ^{\top}\mbox{diag}\{\bsb{y}-\bsbX\bsbb'-\bar\bsbZ \vect(\bsbPhi')\}\bsbZ, \bsbPhi- \bsbPhi' \rangle
+  P_1(\bsbOmega;\bsblambda_{b},\bsbLambda_{\Phi})
+P_2(\bsbOmega;\bsblambda_{\Omega})$. To simplify the notation, we sometimes write them as  $f(\bsb{\Omega})$ and $g(\bsb{\Omega}, \bsb{\Omega}')$ when there is no ambiguity.

Given  $\bsbOmega'$, by simple algebra, the minimization of $g(\bsbOmega,\bsb{\Omega'};\bsblambda_{b},\bsbLambda_{\Phi},\bsblambda_{\Omega})$  subject to $\bsbPhi=\bsbPhi^\top$  reduces to
\begin{align}
\min_{\bsbOmega\in \mathbb{R}^{(p+1) \times p}} \frac{1}{2} \| \bsbOmega - \bsbXi\|_F^2 + P_1(\bsbOmega;\bsblambda_{b}/\tau^2,\bsbLambda_{\Phi}/\tau^2)+P_2(\bsbOmega;\bsblambda_{\Omega}/\tau^2)+ \iota_{\{\bsbPhi=\bsbPhi^{\top}\}}, \label{innLoopObj}
\end{align}
where  $\bsb{\Xi}= [\bsb{\Xi}_{b}, \bsb{\Xi}_{{\Phi}}^{\top}]^{\top}$ with $\bsb{\Xi}_{b}= \bsbb'+\bsbX^{\top}(\bsb{y}-\bsbX\bsbb'- \bar\bsbZ \vect(\bsbPhi') )/\tau^2$ and
$\bsb{\Xi}_{\Phi}= \bsbPhi'+\bsbZ^{\top}\mbox{diag}\{\bsb{y}-\bsbX\bsbb'-\bar\bsbZ\vect(\bsbPhi')\}\bsbZ/\tau^2$, and
$
\iota_{\{\bsbPhi=\bsb{\bsbPhi}^{\top}\}} =
0 $ if $\bsbPhi={\bsbPhi}^{\top}$ and  $+\infty$ otherwise,
which is a lower semicontinuous convex function.

\begin{lemma}
\label{lm_ConvIterProcedure}
Let $\bsb{P} \gets \bsb{0}$, $\bsb{Q}\gets \bsb{0}$, $\bsbXi\gets \bsbXi^0$. Consider the following iterative procedure:
 \begin{equation}
\begin{aligned}
&{\bsbOmega}[:, k]\gets \vec{\bsb{\Theta}}_{S}(\bsbXi[:,k]+\bsb{P}[:, k];\bsb{\lambda}_{\Omega}[k]),  \forall k: 1 \leq k \leq p,\\
&\bsb{P}\gets\bsb{P}+\bsb{\Xi}-\bsb{\Omega},\\ &\bsb{\Xi}[1,:]\gets\Theta_{S}(\bsb{\bsbOmega}[1,:]+\bsb{Q}[1,:];\bsblambda_{b}),\\
&\bsb{\Omega}[2:\mbox{end},:] \gets(\bsb{\Omega}[2:\mbox{end},:] + \bsb{\Omega}[2:\mbox{end},:]^{\top})/2,\\
&\bsb{\Xi}[2:\mbox{end},:]\gets\Theta_{S}(\bsb{\Omega}[2:\mbox{end},:]+\bsb{Q}[2:\mbox{end},:]; (\bsbLambda_{\Phi}+\bsbLambda_{\Phi}^\top)/2),\\
& \bsb{Q} \gets\bsb{Q}+\bsb{\Omega}-\bsb{\Xi}, \label{iterproce}
\end{aligned}
\end{equation}
till convergence. Then both the sequence of      $\bsb{\Xi}$ and the sequence of      $\bsb{\Omega}$
  converge to the globally optimal solution to
$\min_{\bsbOmega \in \mathbb{R}^{(p+1) \times p}}$
$\frac{1}{2}\|\bsbOmega-\bsb{\Xi}^{0}\|_F^2  +\|\bsblambda_{b} \odot \bsbOmega_{b}\|_1+\|\bsbLambda_{\Phi} \odot \bsbOmega_{\Phi}\|_1 +  \|(\bsblambda_{\Omega} \bsb{1}^\top)  \odot \bsbOmega^{\top}\|_{2,1}$ s.t. $\bsbOmega_{\Phi}=\bsbOmega_{\Phi}^{\top}.
$\end{lemma}

Now we prove the
convergence of $\bsbOmega^{(i)}$ in Algorithm \ref{alg1}. First, we notice the following fact.

\begin{lemma}\label{nonexp}
Given any $\lambda\geq 0$, $\vec{\bsb\Theta}_S(\cdot;\lambda)$ is non-expansive, that is,
$\| \vec{\bsb\Theta}(\bsb{a}; \lambda) - \vec{\bsb\Theta}(\tilde {\bsb{a}}; \lambda) \|_2 \leq \|\bsb{a} -
\tilde {\bsb{a}}\|_2$ for all $\bsb{a}, \tilde {\bsb{a}}\in {\mathbb R}^{p}$.
\end{lemma}

Therefore, the mapping of the iteration in Algorithm \ref{alg1} is non-expansive.
 We
use the tool provided by    \cite{opial1967weak} 
for nonexpansive operators to
prove the strict convergence of $\bsbOmega^{(i)}$.
First, the fix point set of the mapping  is non-empty due to the convexity
and the KKT conditions. The mapping  is
also {asymptotically regular in the sense that}  for any starting point $\bsbOmega^{(0)}$, $\| \bsbOmega^{(i+1)} -\bsbOmega^{(i)} \| \rightarrow 0$
as $i\rightarrow \infty$.
This property is implied by the following lemma with   $\tau>\| \breve \bsbZ\|_2$.

\begin{lemma}
\label{lm_FuncValueDecreIter}
For the sequence $\{\bsb{\Omega}^{(i)}\}$  generated by Algorithm \ref{alg1},
\begin{equation}
\begin{aligned}
f(\bsb{\Omega}^{(i)};\bsblambda_{b},\bsbLambda_{\Phi},\bsblambda_{\Omega})-f(\bsb{\Omega}^{(i+1)};\bsblambda_{b},\bsbLambda_{\Phi},\bsblambda_{\Omega})
&\geq \frac{1}{2}(\tau^2-\| \breve{\bsbZ}\|_2^2)\|\bsbOmega^{(i+1)}-\bsbOmega^{(i)}  \|_F^2,\label{L11}
\end{aligned}
\end{equation}
for any $i\geq 0$. Moreover,  $\{\bsbOmega^{(i)}\}$ is uniformly bounded.
\end{lemma}

 \eqref{L11} actually does not require   the optimality of $\bsb{\Omega}^{(i+1)}$  in  minimizing  $g(\cdot, \bsb{\Omega}^{(i)})$ (cf.  \eqref{innLoopObj}). From the proof of Lemma \ref{lm_FuncValueDecreIter}, any $\bsb{\Omega}^{(i+1)}$  satisfying $g(\bsbOmega^{(i+1)}, \bsb{\Omega}^{(i)}) \le g(\bsbOmega^{(i)},\bsb{\Omega}^{(i)})(=f(\bsbOmega^{(i)}))$ makes \eqref{L11} hold. There is no need  to run the inner loop till convergence. Moreover,  to get the asymptotic regularity and the uniform boundedness of   $\{\bsbOmega^{(i)}\}$, it  suffices to  terminate the inner loop when  $g(\bsbOmega^{(i+1)}, \bsb{\Omega}^{(i)}) \le g(\bsbOmega^{(i)},\bsb{\Omega}^{(i)}) +  \theta \|\bsbOmega^{(i+1)}-\bsbOmega^{(i)}  \|_F^2/2
$ or
$$
\frac{g(\bsbOmega^{(i+1)}, \bsb{\Omega}^{(i)}) - g(\bsbOmega^{(i)},\bsb{\Omega}^{(i)})}{\|\bsbOmega^{(i+1)}-\bsbOmega^{(i)}  \|_F^2/2} \le   \theta,
$$
where $\theta$ is a pre-specified constant satisfying  $-\tau^2\le \theta< \tau^2-\| \breve{\bsbZ}\|_2^2$. When  $\bsbOmega^{(i+1)}$ exactly solves  \eqref{innLoopObj}, $\theta=-\tau^2$ is satisfied.

With all  Opial's conditions  satisfied, the sequence $\{\bsbOmega^{(i)}\}$ has a unique limit point $\bsbOmega^* = [\bsbb^*,\bsbPhi^{* \top}]^{\top}$, and $\bsbOmega^*$ is a fixed point of Algorithm \ref{alg1}.
\

Next we prove the optimality of any fixed point $\bsbOmega^*$ regarding  problem \eqref{general0}. By Lemma \ref{lm_ConvIterProcedure}, $\bsbOmega^*$ satisfies the Karush-Kuhn-Tucker (KKT) conditions of \eqref{innLoopObj}:
 \begin{equation*}
 \begin{cases}
 \boldsymbol{0} \in \bsbOmega^* - \bsb{\Xi}  + \partial P_1(\bsb{\Omega};\bsblambda_{b}/\tau^2,\bsbLambda_{\Phi}/\tau^2)+ \partial P_2(\bsb{\Omega};\bsbLambda_{\Omega}/\tau^2) +  \partial h(\bsbPhi^*), \\
 \bsb{\bsbPhi}^*- \bsb{\bsbPhi}^{*\top}  = \bsb{0},
 \end{cases}
 \end{equation*}
where  $h(\bsbPhi;
\bsb{L})=\langle \bsb{L}, \bsbPhi - \bsbPhi^{\top}\rangle$.
With the fixed point property established above, substituting $\bsb{\Xi}$ with $\bsb{\Xi}_{b}=\bsbb^{*}+\bsbX^{\top}(\bsb{y}-\bsbX\bsbb^{*}-\bar{\bsbZ}\vect(\bsbPhi^*))/\tau^2$ and
 $\bsb{\Xi}_{\Phi}=\bsbPhi^{*}+\bsbZ^{\top}\mbox{diag}\{\bsb{y}-\bsbX\bsbb^{*}-\bar{\bsbZ}\vect(\bsbPhi^*)\}\bsbZ/\tau^2$, we have
\begin{align*}
\begin{cases}
\bsb{0} &\in  \bsbX^{\top}(\bsb{y}-\bsbX\bsbb^{*}-\bar{\bsbZ}\vect(\bsbPhi^*))
     + \partial_{\bsbb} P_1(\bsb{\Omega};\bsblambda_{b},\bsbLambda_{\Phi})+ \partial_{\bsbb} P_2(\bsb{\Omega};\bsbLambda_{\Omega}),\\
\bsb{0} &\in \bsbZ^{\top}\mbox{diag}\{\bsb{y}-\bsbX\bsbb^{*}-\bar{\bsbZ}\vect(\bsbPhi^*)\}\bsbZ+ \partial_{\bsbPhi} (P_1(\bsb{\Omega};\bsblambda_{b},\bsbLambda_{\Phi})+P_2(\bsb{\Omega};\bsbLambda_{\Omega}))\\
&+ \partial h(\bsbPhi^*; \tau^2 \bsb{L}),\\
\bsb{0} &=\bsbPhi^*-\bsbPhi^{*\top},
\end{cases}
\end{align*}
which are exactly the KKT conditions for the convex problem  \eqref{general0}. Hence $\bsbOmega^*$ is a global minimizer of problem  \eqref{general0}.

Similarly, the global optimality of $\lim_{i\rightarrow\infty} \bsbOmega^{(i)}=\bsbOmega^*$ does not require solving   \eqref{innLoopObj} exactly.
In the following, we use $\|\bsb{A} \|_{\bsb{\Sigma}}$ to denote $\{tr(\bsb{A}^\top \bsb{\Sigma} \bsb{A})\}^{1/2}$. Let $\bsbOmega^o$ be an arbitrary optimal solution to  \eqref{general0}.
From the proof of Lemma \ref{lm_FuncValueDecreIter}, the optimal $\bsbOmega^{(i+1)}$ satisfies  $g(\bsbOmega^{(i+1)}, \bsb{\Omega}^{(i)}) + \frac{\tau^2}{2} \| \bsbOmega^{(i+1)} -  \bsb{\Omega}^o\|_F^2 \leq g(\bsbOmega^o,\bsb{\Omega}^{(i)})$.
We can  relax it to
$$
g(\bsbOmega^{(i+1)}, \bsb{\Omega}^{(i)}) + \frac{1}{2} \| \bsbOmega^{(i+1)} -  \bsb{\Omega}^o\|_{\tau^2 \bsb{I}- \theta \bsbX^\top \bsbX}^2 \leq g(\bsbOmega^o,\bsb{\Omega}^{(i)})
$$
for some $\theta: 0\le \theta\le 1$.
Then by the construction of $g$, we obtain
\begin{align*}
f(\bsbOmega^{(i+1)})  - f(\bsb{\Omega}^{o}) &\le\frac{1}{2} \| \bsbOmega^{(i)} -  \bsb{\Omega}^o\|_{\tau^2 \bsb{I} -  \bsbX^\top \bsbX}^2  -  \frac{1}{2} \| \bsbOmega^{(i+1)} -  \bsb{\Omega}^o\|_{\tau^2 \bsb{I}- \theta \bsbX^\top \bsbX}^2 \\& \quad-   \frac{1}{2} \| \bsbOmega^{(i+1)} -  \bsb{\Omega}^{(i)}\|_{\tau^2 \bsb{I}-  \bsbX^\top \bsbX}^2 \\
& \le \| \bsbOmega^{(i)} -  \bsb{\Omega}^o\|_{\tau^2\bsb{I} -  \bsbX^\top \bsbX}^2    -  \frac{1}{2} \| \bsbOmega^{(i+1)} -  \bsb{\Omega}^o\|_{\tau^2\bsb{I} - \theta \bsbX^\top \bsbX}
\\
&\le \frac{1}{2} \| \bsbOmega^{(i)} -  \bsb{\Omega}^o\|_{\tau^2 \bsb{I}- \theta \bsbX\top \bsbX}^2  -  \frac{1}{2} \| \bsbOmega^{(i+1)} -  \bsb{\Omega}^o\|_{\tau^2\bsb{I} - \theta \bsbX\top \bsbX}.
\end{align*}
It follows that
\begin{align}
{\bar f}_{i+1}
 - f(\bsb{\Omega}^{o}) \le \frac{1}{2(i+1)}  \| \bsbOmega^{(0)} -  \bsb{\Omega}^o\|_{\tau^2 \bsb{I}- \theta \bsbX\top \bsbX}^2,
\end{align}
where $\bar f_{i+1} = \frac{1}{i+1}\sum_{k=1}^{i+1} f(\bsbOmega^{(k)})$.
Letting  $i\rightarrow \infty$, we get $f(\bsbOmega^*)=f(\lim_{i\rightarrow \infty} \bsbOmega^{(i)})=\lim_{i\rightarrow \infty} {\bar f}_{i}=f(\bsbOmega^o)$, and so $\bsbOmega^*$ must also be a global minimizer of     \eqref{general0}.
\paragraph*{Proof of Lemma \ref{lm_ConvIterProcedure}}
First, we consider
$$
\min_{\bsbOmega=[\bsbb, \bsbPhi]^\top}\frac{1}{2}\|\bsbOmega-\bsb{\Xi}^{0}\|_F^2  +\|\bsblambda_{b} \odot \bsbb\|_1+\|\bsbLambda_{\Phi} \odot \bsb\Phi\|_1 +\iota_{\{\bsbPhi=\bsb{\bsbPhi}^{\top}\}}. $$
Let $\bsbPhi^0 = \bsbXi_{\Phi}^0$, $\bsbb^0 = \bsbXi_b^0$. By some simple algebra,
$
\|\bsbPhi - \bsbPhi^0\|_F^2 = \|\bsbPhi - (\bsbPhi^0+\bsbPhi^{0\top})/{2}\|_F^2 + C(\bsbPhi^0),
$
where the last term does not depend on $\bsbPhi$. Based on Lemma 1 of  \cite{Group_She}, it is not difficult to show that the  global minimizer is given by $\bsbb_o=\Theta_{S}(\bsbb^{0};\bsblambda_{b})$,   $\bsbPhi_o=\Theta_{S}( (\bsbPhi^0+\bsbPhi^{0\top})/{2};(\bsbLambda_{\Phi}+\bsbLambda^{\top}_{\Phi})/{2})$.

Similarly, the globally optimal  solution to $\min_{\bsbOmega=[\bsbb, \bsbPhi]^\top}\frac{1}{2}\|\bsbOmega-\bsb{\Xi}^{0}\|_F^2  +\|(\bsblambda_{\Omega} \bsb{1}^\top)  \odot \bsbOmega^{\top}\|_{2,1} $ is given by   $ \bsbOmega_o$ with  $ \bsbOmega_o[;,k] =\vec{\bsb{\Theta}}_S(\bsbXi^{0}[:,k];\bsblambda_{\Omega}[k] )$, $1\leq k\leq p$.
Applying  Theorem 3.2 and 3.3 in \cite{Bauschke} guarantees the strict convergence
of the iterates and the global optimality of the limit point.   \qed

\paragraph*{Proof of Lemma \ref{nonexp}} The conclusion follows from the fact that $\vec {\bsb{\Theta}}$ is a proximity operator associated with a convex function.
%
%
\qed

\paragraph*{Proof of Lemma \ref{lm_FuncValueDecreIter}}
From the optimality of $\bsbOmega^{(i+1)}$ and the convexity of the penalties in  \eqref{innLoopObj},  it is easy to see that for any $\bsbOmega$,
\begin{align}
g(\bsbOmega^{(i+1)}, \bsb{\Omega}^{(i)};  \bsblambda_{b},\bsbLambda_{\Phi},\bsbLambda_{\Omega}) + \frac{\tau^2}{2} \| \bsbOmega^{(i+1)} -  \bsb{\Omega}\|_F^2 \leq g(\bsbOmega,\bsb{\Omega}^{(i)};  \bsblambda_{b}, \bsbLambda_{\Phi},\bsblambda_{\Omega}).
\end{align}
It follows from  the construction of $g$ that
$$
 g(\bsbOmega^{(i+1)}, \bsb{\Omega}^{(i)};  \bsblambda_{b},\bsbLambda_{\Phi},\bsbLambda_{\Omega}) \leq g(\bsbOmega^{(i)},\bsb{\Omega}^{(i)};  \bsblambda_{b}, \bsbLambda_{\Phi},\bsblambda_{\Omega}) =f(\bsb{\Omega}^{(i)};  \bsblambda_{b}, \bsbLambda_{\Phi},\bsblambda_{\Omega}).
$$
On the other hand, noticing that the gradient of $l$ with respect to $\bsbPhi$ is $\bsbZ^{\top}\mbox{diag}\{\bsb{y}-\bsbX\bsbb'-\bar\bsbZ \vect(\bsbPhi'))\}\bsbZ$, Taylor expansion yields $$g(\bsbOmega, \bsb{\Omega}';  \bsblambda_{b},\bsbLambda_{\Phi},\bsbLambda_{\Omega})\ge f(\bsb{\Omega};  \bsblambda_{b}, \bsbLambda_{\Phi},\bsblambda_{\Omega})+ \frac{(\tau^2-\| \breve{\bsbZ}\|_2^2)}{2}\|\bsbOmega'-\bsbOmega\|_F^2$$ for any $\bsbOmega, \bsbOmega'$.

Furthermore,
$
\| (\bsblambda_{\Omega}\bsb{1}^\top) \odot\bsb{\Omega}^{(i)\top} \|_{2,1} \leq f(\bsbOmega^{(i)};\bsblambda_{b},\bsbLambda_{\Phi}) \leq f(\bsbOmega^{(0)};\bsblambda_{b},\bsbLambda_{\Phi}).
$ Because $\bsblambda_{\Omega}>\bsb{0}$, the conclusion follows. \qed

\subsection{Proof of Theorem \ref{th_oracle}}
\label{sec:appOracle}
In this proof,  we use  $C$, $c$ to denote universal constants. They are not necessarily the same at each occurrence. Throughout the proof, $P(\bsb{b},\bsbPhi;\lambda_1,\lambda_2)$ is short for $\lambda_1\|\bsbPhi\|_1+\lambda_2\|\bsbOmega^{\top}\|_{2,1}$, $\lambda_1'$ is short for $\lambda_1\|\breve{\bsbX}\|_2$, and $\lambda_2'$ is short for $\lambda_2\|\breve{\bsbX}\|_2$.
We prove the theorem under a less restrictive condition.

First, because  $\bsb{\hat{\Omega}}$ is a global minimizer of \eqref{ob1} and $P$ is convex,  for any  $(\bsbb, \bsbPhi) \in \mathbb{R}^p \times  \mathbb{R}^{p \times p}$ with $\bsbPhi$ symmetric (not necessarily satisfying SH),  we have
$$
\ell(\hat{\bsbb}, \hat{\bsbPhi} )+P(\hat{\bsbb},\hat{\bsbPhi};\lambda_1',\lambda_2') +  \frac{1}{2} M(\bsb{{b}}-\hat\bsbb, \bsb{{\bsbPhi}}-\hat\bsbPhi)\leq \ell(\bsbb,\bsbPhi)+P(\bsbb,\bsbPhi;\lambda_1',\lambda_2').
$$
 Based on the model assumption,  it is easy to see that
\begin{align*}
&\frac{1}{2}M(\bsbb^*-\hat{\bsbb},\bsbPhi^*-\hat{\bsbPhi})+P(\hat{\bsbb},\hat{\bsbPhi};\lambda_1',\lambda_2')+ \frac{1}{2} M(\bsbDelta_{b}, \bsbDelta_{\Phi}) \\\leq &  \frac{1}{2}M(\bsbb^*-\bsbb,\bsbPhi^*-\bsbPhi)+P(\bsbb,\bsbPhi;\lambda_1',\lambda_2')
+\langle \bsb{\varepsilon},\bsb{X}\bsbDelta_{b}+\mbox{diag}(\bsb{X}\bsbDelta_{\Phi}\bsb{X}^{\top}) \rangle,
\end{align*}
 where $\bsbDelta_{b}=\hat{\bsb{{b}}}-\bsbb$, $\bsbDelta_{\Phi}=\hat{{\bsbPhi}}-\bsbPhi$.

Introduce $\breve{\bsbDelta}=[\bsbDelta_{b}^{\top}, (\vect(\bsbDelta_{\Phi}))^{\top}]^{\top}$ and  $\tilde{\mathcal{J}}=\mathcal{J} \cup \hat{\mathcal{J}}$, where   $\mathcal{J}:=\{j:\vect({\bsbOmega})_j \neq 0\}$.  and     $\hat{\mathcal{J}}:=\{j:\vect{(\hat{\bsbOmega})}_j \neq 0\}$. Denote by $\Proj_{\tilde{\mathcal J}}$ the orthogonal projection matrix onto  the column space of  ${{\breve \bsbX}_{\tilde{\mathcal J}}}$.
Then
\begin{align*}
&\langle \bsbeps, \bsbX\bsbDelta_{b}+\mdiag(\bsbX\bsbDelta_{\Phi}\bsbX^{\top})\rangle \\
=& \langle \bsbeps,\breve{\bsbX}\breve{\bsbDelta}\rangle \\
=&\langle \bsbeps, \Proj_{\tilde{\mathcal{J}}}\breve{\bsbX}_{\tilde{\mathcal{J}}}\breve{\bsbDelta}_{\tilde{\mathcal{J}}}\rangle \\
=& \langle \Proj_{\tilde{\mathcal{J}}}\bsbeps,\bsbX\bsbDelta_{b}+\mdiag(\bsbX\bsbDelta_{\Phi}\bsbX^{\top}) \rangle\\
\leq & \frac{1}{2a_1} M(\bsbDelta_{b},\bsbDelta_{\Phi})+\frac{a_1}{2}\|\Proj_{\tilde{\mathcal{J}}}\bsbeps\|_2^2,
\end{align*}   
for any $a_1>0$. Representing   $\|\hat \bsbOmega\|_{2,1}$ as $\|(\hat \bsbOmega_{\mathcal J_G})^{\top}\|_{2,1}+\|(\hat \bsbOmega_{\mathcal J_G^c})^{\top}\|_{2,1}$ and  $\|\hat \bsbPhi\|_1$ as $\|\hat \bsbPhi_{\mathcal J_e}\|_1+\|\hat \bsbPhi_{\mathcal J_e^c}\|_1$,   respectively, and applying the triangle inequalities  of $\| \cdot \|_1$ and $\|\cdot\|_{2,1}$, we have
\begin{equation}\begin{aligned}
&\frac{1}{2} M(\bsbb^*-\hat{\bsbb},\bsbPhi^*-\hat{\bsbPhi})+\lambda_1'\|(\bsbDelta_{\Phi})_{\mathcal{J}^c_e}\|_1+\lambda_2'\|(\bsbDelta_{\mathcal{J}^c_G})^{\top}\|_{2,1}+ \frac{1}{2} M(\bsbDelta_{b}, \bsbDelta_{\Phi})\\
\leq &\frac{1}{2} M(\bsbb^*-\bsbb,\bsbPhi^*-\bsbPhi)+\frac{1}{2a_1} M(\bsbDelta_{b},\bsbDelta_{\Phi})+\lambda_1'\|(\bsbDelta_{\Phi})_{\mathcal{J}_e}\|_1+ \lambda_2'\|(\bsbDelta_{\mathcal{J}_G})^{\top}\|_{2,1}
\\ &+\frac{a_1}{2}\|\Proj_{\tilde{\mathcal{J}}}\bsbeps\|_2^2.
\end{aligned}\label{basicineq1}
\end{equation}

How to treat the key term  $\|\Proj_{\tilde{\mathcal{J}}}\bsbeps\|_2^2$ is nontrivial. We derive two inequalities from statistical analysis and computational analysis.

First, we bound  $\|\Proj_{\tilde{\mathcal{J}}}\bsbeps\|_2^2$ with some cardinality measures.
The goal of this step  is to show for a large constant $L$,  $\sup_{\mathcal J\subset[p(p+1)]}\{ \|\Proj_{\mathcal J} \bsbeps\|_2^2 -   {L  P_o(J_G(\mathcal J), J_e(\mathcal J))}\}$ is negative with high probability, where $P_o(J_{G}, J_{e}) =\sigma^2 \{J_G \log (ep/J_G)+J_e\log(e p J_G / J_e)\}$. (If the reference signal $(\bsb{b}, \bsb{\Phi})$  satisfies SH, so does $\tilde {\mathcal J}$, for which  $P_o$ can be strengthened to $\sigma^2 \{J_G \log (ep/J_G)+J_e\log(e J_G^2 / J_e)\}$ following the lines of the proof of Lemma \ref{noisebound}.)

For notational convenience, we extend the \textit{J}-measures defined in \eqref{Jdefs} to an index set: Given any  $\mathcal J\subset [p(p+1)]$, with  $\bsbOmega=[\bsbb, \bsbPhi^\top]^\top$   a binary matrix satisfying $(\vect(\bsbOmega))_{\mathcal J}=1$ and $(\vect(\bsbOmega))_{\mathcal J^c}=0$,   define $\mathcal J_G (\mathcal J)=\mathcal J_{G}(\bsbOmega)$, $\mathcal J_e(\mathcal J) = \mathcal J_e(\bsbPhi)$.
Given $\mathcal J\subset [p(p+1)]$, $\mathcal J_G(\mathcal J)\subset [p]$ and $\mathcal J_e(\mathcal J)\subset [p^2]$.

\begin{lemma}\label{noisebound}
Given $1\leq J_G \leq  p$ and $0\leq J_e  \leq p^2$, for  any $t\geq 0$,
\begin{align}
\EP \left(\sup_{|\mathcal J_G|=J_G, |\mathcal J_e|=J_e} \|\Proj_{\mathcal J} \bsbeps\|_2^2\geq    {L_0  P_o(J_G, J_e)}+ t \sigma^2 \right) \leq C\exp(- c t),
\end{align}
where $L_0, C, c>0$ are universal constants.
Here, we omit the dependence of $\mathcal J_G$ and $\mathcal J_e$ on $\mathcal J$  for brevity.\end{lemma}

Define
\[
R_{J_G, J_e}= \left( \sup_{|\mathcal J_G|=J_G, |\mathcal J_e|=J_e}\|\Proj_{\mathcal J} \bsbeps\|_2^2 -   {L \sigma^2 \{J_e\log(ep)+J_G\log(ep)\} }\right)_{+},
  \]
and  $R=\sup_{1\leq J_G \leq  p, 0\leq J_e  \leq p^2} R_{J_G, J_e}$ (when $J_G(\mathcal J)=0$, $\Proj_{\mathcal J} \bsbeps=\bsb{0}$).   Then \begin{align}
\|\Proj_{\tilde{\mathcal{J}}}\bsbeps\|_2^2
\leq& L \sigma^2 \{{\tilde J}_e \log (e p) +{\tilde J}_G \log (e p )\} +R \notag \\
\leq & L \sigma^2 (J_e +  J_G) \log (ep) +  L \sigma^2 \{{\hat J}_e \log (e p) +{\hat J}_G \log (e p )\} +R .  \label{seesaw1}
\end{align}
   Set  $L> L_0$.  Then, by Lemma \ref{noisebound} and $J_e\ge J_G$,
\begin{align*}
& \EP(R\geq t \sigma ) \\
\leq & \sum_{J_G=1}^p \sum_{J_e=0}^{p^2} \EP(R_{J_G,J_e}\geq t\sigma )\\
 \leq &  \sum_{J_G=1}^p \sum_{J_e=0}^{p^2}  C \exp(-c t) \exp\left\{- c \left (L-L_0) (J_e\log(ep)+J_G\log(ep)\right)\right\} \\
 \leq & C' \exp(-c t^2),
\end{align*}
from which it follows that  $\EE R \leq C \sigma^2$.
It is also easy to see that for sufficiently large $L$,  $R\leq 0$ occurs with probability at least $1-Cp^{-c L}$.

Next, we derive an inequality based on  the computational optimality of $(\hat\bsbb, \hat \bsbPhi)$.  Due to the convexity of the problem, $(\hat\bsbb, \hat \bsbPhi)$ is a  stationary point of
\begin{align*}
&\frac{1}{2}\|\bsb{y}-\bsbX\bsbb-\mbox{diag}(\bsbX\bsbPhi\bsbX^{\top})\|_2^2 +\lambda_1'\| \bsbPhi\|_1 \\
& \quad + \lambda_2' \| \bsbOmega^{\top} \|_{2,1} + \sum_{j<k} l_{j,k} (\phi_{j,k} -\phi_{k,j}), \end{align*}
where $l_{j,k}$ are Lagrangian multipliers.
By the  KKT conditions, for any $j\in \hat {\mathcal J}_G$,  $\hat{b}_j$ is nonzero  (with probability 1) and so satisfies \begin{equation}
\lambda_2'\frac{\hat{b}_j}{\|\hat{\bsb{\bsbOmega}}_j\|_2}= -\bsbx_j^{\top}(\breve{\bsbX} \vect(\hat{\bsb{\bsbOmega}})-\bsb{y}). \label{kkt_b}
\end{equation}
Similarly, for any $j'\in \hat{\mathcal J}_e$, letting $j_c = \lceil \frac{j'}{p}\rceil$, $j_r=j'-(j_c -1)p$,  $\hat{\phi}_{j_r, j_c}$ satisfies  $\lambda_1'\msgn(\hat{\phi}_{j_r, j_c})+ \lambda_2'\frac{\hat{\phi}_{j_r, j_c}}{\|\hat{\bsb{\bsbOmega}}_{j_c}\|_2} + l_{j_r, j_c} \msgn(j_r-j_c)=- \bar{\bsbX}_{j'}^{\top}(\breve{\bsbX} \vect(\hat{\bsb{\bsbOmega}})-\bsb{y})$, where $\msgn(t)=1, 0, -1$ for $t>0$, $=0$,  $<0$, respectively. Adding the equations for $\hat{\phi}_{j_r, j_c}$ and $\hat{\phi}_{j_c, j_r}$ can cancel $l_{j_r,j_c}$. So when $j_r\neq j_c$, we have
\begin{align}
2\lambda_1'\msgn(\hat{\phi}_{j_r, j_c})+ \lambda_2'\frac{\hat{\phi}_{j_r, j_c}}{\|\hat{\bsb{\bsbOmega}}_{j_c}\|_2} +\lambda_2'\frac{\hat{\phi}_{j_r, j_c}}{\|\hat{\bsb{\bsbOmega}}_{j_r}\|_2} = -2\bar{\bsbX}_{j'}^{\top}(\breve{\bsbX} \vect(\hat{\bsb{\bsbOmega}})-\bsb{y}), \label{kkt_phi_offdiag}
\end{align}
and when $j_c=j_r$, we have
\begin{align}
\lambda_1'\msgn(\hat{\phi}_{j_r, j_c})+ \lambda_2'\frac{\hat{\phi}_{j_r, j_c}}{\|\hat{\bsb{\bsbOmega}}_{j_c}\|_2} = -\bar{\bsbX}_{j'}^{\top}(\breve{\bsbX} \vect(\hat{\bsb{\bsbOmega}})-\bsb{y}). \label{kkt_phi_diag}
\end{align}
Notice that $\msgn(\hat{\phi}_{j_r, j_c})   \hat{\phi}_{j_r, j_c}$ is always non-negative. Squaring both sides of \eqref{kkt_b},  \eqref{kkt_phi_offdiag} and \eqref{kkt_phi_diag} and summing over all $j\in \hat {\mathcal J}_G$ and  $j'\in \hat{\mathcal J}_e$, we obtain
\begin{align*}
\lambda_1'^2  {\hat{J}_{e}}+\lambda_2'^2  \sum_{j\in \hat{\mathcal  J_G}}\frac{{\hat b}_j^2+\sum_{k:\hat\phi_{k,j}\neq 0}\hat\phi_{k, j}^2}{\|\hat{\bsb{\bsbOmega}}_j\|_2^2}  \leq 2\|\breve{\bsbX}_{\hat{\mathcal J}}^{\top}(\breve{\bsbX} \vect(\hat{\bsb{\bsbOmega}})-\bsb{y})\|_2^2,
\end{align*}
and so
$\lambda_1'^2  {\hat{J}_{e}}+\lambda_2'^2    \hat{J}_G\leq 2\|\breve{\bsbX}_{\hat{\mathcal J}}^{\top}(\breve{\bsbX} \vect(\hat{\bsb{\bsbOmega}})-\bsb{y})\|_2^2$. It follows from  $\hat{\mathcal J} \subset \tilde{\mathcal J}$ and the Cauchy-Schwarz inequality  that
\begin{align}
\lambda_1^2\hat{J}_{e}+\lambda_2^2\hat{J}_{G} &
 \leq 3 M(\bsbb^*-\hat{\bsbb},\bsbPhi^*-\hat{\bsbPhi})+6\|\Proj_{\tilde{\mathcal{J}}}\bsbeps\|_2^2.\label{seesaw2}
\end{align}

Now, combining the optimization  inequality  \eqref{seesaw2} and the statistical inequality \eqref{seesaw1} can yield a bound of $\|\Proj_{\tilde{\mathcal{J}}}\bsbeps\|_2^2$. In fact, with $\lambda^o=\sigma\sqrt{\log (e p)}$, $\lambda_1=  = A_1\lambda^o$, $\lambda_2 = A_2\lambda^o$, $  A_i\ge A$,  $A^2= A_0 L$ for some $A_0>0$,
we get
\begin{align*}
(1-\frac{6}{A_0})\|\Proj_{\tilde{\mathcal{J}}}\bsbeps\|_2^2 \leq L\sigma^2 (J_e + J_G) \log (ep) + \frac{3}{A_0}  M(\bsbb^*-\hat{\bsbb},\bsbPhi^*-\hat{\bsbPhi})+R.
\end{align*}
Plugging this into \eqref{basicineq1}
gives \begin{align*}
&\frac{1}{2} M(\bsbb^*-\hat{\bsbb},\bsbPhi^*-\hat{\bsbPhi})+A \lambda^o \|\breve{\bsbX}\|_2\|(\bsbDelta_{\Phi})_{\mathcal{J}^c_e}\|_1+A\lambda^o\|\breve{\bsbX}\|_2\|(\bsbDelta_{\mathcal{J}^c_G})^{\top}\|_{2,1}\\
\leq &\frac{1}{2} \left(1+\frac{3a_{1}}{A_0-6}\right) M(\bsbb^*-\bsbb,\bsbPhi^*-\bsbPhi)+\frac{1}{2} \left(\frac{1}{a_1}-1\right) M(\bsbDelta_{b},\bsbDelta_{\Phi})\\&+A \lambda^o \|\breve{\bsbX}\|_2\|(\bsbDelta_{\Phi})_{\mathcal{J}_e}\|_1+ A \lambda^o \|\breve{\bsbX}\|_2\|(\bsbDelta_{\mathcal{J}_G})^{\top}\|_{2,1}
+\frac{a_1 L}{2(1-6/A_0)}\sigma^2 (J_e + J_G) \log (ep) \\&+ \frac{a_1}{2(1-6/A_0)}R.
\end{align*}
Assume the following condition holds\begin{equation}
\begin{aligned}
 &\|\breve{\bsbX}\|_2(\|(\bsbDelta_{\Phi})_{\mathcal{J}_e}\|_1+\|(\bsbDelta_{\mathcal{J}_G})^{\top}\|_{2,1}) \\\le & \|\breve{\bsbX}\|_2(\|(\bsbDelta_{\Phi})_{\mathcal{J}_e^c}  \|_1+\|(\bsbDelta_{\mathcal{J}_G^c})^{\top}\|_{2,1} )+ K \sqrt{J_e +J_G}\|\bsbX\bsbDelta_{b}+\bar{\bsbX}\vect(\bsbDelta_{\Phi})\|_2,
\end{aligned}\label{comparegcond}
\end{equation}
for some $K$ large enough. Then we  have
\begin{align*}
&\frac{1}{2} M(\bsbb^*-\hat{\bsbb},\bsbPhi^*-\hat{\bsbPhi})\\\leq &\frac{1}{2} \left(1+\frac{3a_{1}}{A_0-6}\right) M(\bsbb^*-\bsbb,\bsbPhi^*-\bsbPhi)+\frac{1}{2} \left(\frac{1}{a_1}-1\right) M(\bsbDelta_{b},\bsbDelta_{\Phi})+ \frac{a_1 R}{2(1-6/A_0)}\\&+ A \lambda^o K  \sqrt{J_e +J_G}\|\bsbX\bsbDelta_{b}+\bar{\bsbX}\vect(\bsbDelta_{\Phi})\|_2+\frac{a_1 L}{2(1-6/A_0)}\sigma^2 (J_e + J_G) \log (ep) \\
\le & \frac{1}{2} \left(1+\frac{3a_{1}}{A_0-6}\right) M(\bsbb^*-\bsbb,\bsbPhi^*-\bsbPhi)+\frac{1}{2} \left(\frac{1}{a_1}+\frac{1}{a_2}-1\right) M(\bsbDelta_{b},\bsbDelta_{\Phi})\\&+\left\{\frac{a_1 }{2(1-6/A_0)}+\frac{a_2 A_0 }{2} K^2\right\}L\sigma^2 (J_e + J_G) \log (ep) + \frac{a_1}{2(1-6/A_0)}R.
\end{align*}
Hence we get $\EE M(\bsbb^*-\hat{\bsbb},\bsbPhi^*-\hat{\bsbPhi})\lesssim M(\bsbb^*-\bsbb,\bsbPhi^*-\bsbPhi) + (1\vee K^2) \sigma^2 (J_e + J_G) \sigma^2 \log (e p)+\sigma^2 $ if we choose the constants $a_1$, $a_2$, $A_0$ satisfying, say, $a_1=2$, $a_2=2$,  $A_0\ge 7$.

To complete the proof, we show that $\mathcal{A}(\mathcal{J}_e,\mathcal{J}_G,\vartheta,\kappa)$ with  $\kappa>0$ and   $\vartheta\ge0$ implies  \eqref{comparegcond}.  Consider two cases. The case $(1+\vartheta)(\|(\bsbDelta_{\Phi})_{\mathcal{J}_e}\|_1+\|(\bsbDelta_{\mathcal{J}_G})^{\top}\|_{2,1})\le \|(\bsbDelta_{\Phi})_{\mathcal{J}_e^c}  \|_1+\|(\bsbDelta_{\mathcal{J}_G^c})^{\top}\|_{2,1}   $ is trivial. Suppose the reverse inequality holds. Then
\begin{align*}
& \|\breve{\bsbX}\|_2 (\|(\bsbDelta_{\Phi})_{\mathcal{J}_e}\|_1+\|(\bsbDelta_{\mathcal{J}_G})^{\top}\|_{2,1}) \\\le &  \|\breve{\bsbX}\|_2(J_e + J_G)^{1/2} ( \|(\bsbDelta_{\Phi})_{\mathcal{J}_e}  \|_2^2 +  \|(\bsbDelta_{\mathcal{J}_G})^{\top}\|_{F}^2    )^{1/2} \\
  \le &  (J_e + J_G)^{1/2} \| \bsbX\bsbDelta\|_F/\kappa^{1/2},
\end{align*}
and so  \eqref{comparegcond} holds with  $K=1/\kappa^{1/2}$.

Finally, we can also plug \eqref{seesaw1} into \eqref{seesaw2}, resulting in
\begin{equation}
\lambda_1^2 \mathbb{E}[\hat{J}_e]+ \lambda_2^2 \mathbb{E}[\hat{J}_G] \lesssim M(\bsbb^*-\bsbb,\bsbPhi^*-\bsbPhi)+ \lambda_1^2J_e+\lambda_2^2 J_G+\sigma^2,
\end{equation}
and so \eqref{cardbound} follows.

The proof for Theorem 2' follows the same lines; the details are omitted.\\

\paragraph*{Proof of Lemma \ref{noisebound}}
First, notice that for fixed  ${\mathcal J}$, $\|\Proj_{{\mathcal{J}}}\bsbeps\|_2^2 /\sigma^2\sim \chi^2(D)$ with $D\leq J$. The  standard tail bound for $X\sim \chi^2(D)$ gives  $\EP(X - D \geq t)\leq \exp(-{t^2}/(4{(D + t)}))$,  $\forall t\geq 0$  \cite[Lemma 1]{laurent2000adaptive}.

It is easy to see that
\[\sup_{|\mathcal J_G|=J_G, |\mathcal J_e|=J_e} \|\Proj_{\mathcal J} \bsbeps\|_2^2\leq
\sup_{J^{11}+J^{10}+J^{01}=J_G}\sup_{  {\mathcal Q}(J^{11},J^{10},J^{01}) }  \sup_{|\mathcal{J}_e|=J_e, \mathcal{J}_e\subset \mathcal J_G^{(e)}}\|\Proj_{\mathcal J} \bsbeps\|_2^2,
\] where   $\mathcal{Q}(J^{11},J^{10},J^{01}) = \{(\mathcal{J}^{11},\mathcal{J}^{10},\mathcal{J}^{01} )\subset [p]^3: |\mathcal{J}^{11}|=J^{11}, |\mathcal{J}^{10}|=J^{10}, |\mathcal{J}^{01}| =J^{01}, \mathcal{J}^{11}\cap\mathcal{J}^{10}=\mathcal{J}^{10}\cap\mathcal{J}^{01}=\mathcal{J}^{11}\cap\mathcal{J}^{01}=\emptyset \}$ and $\mathcal J_G^{(e)}$ denotes the index set $\{1\leq j \leq p^2: \lceil \frac{j}{p}\rceil  \in {\mathcal J}_G  \}$.  ${\mathcal Q}(J^{11},J^{10},J^{01})$ contains  $${p \choose J^{11},\ J^{10}, \ J^{01}, \ p-J^{11}-J^{10}-J^{01}}$$  many elements. This multinomial coefficient is bounded by $p \choose J_G$, and by Stirling's approximation, $\log {p \choose J_G} \leq J_G \log (ep/J_G)$. Similarly,
  $\sup_{|\mathcal{J}_e|=J_e, \mathcal{J}_e\subset \mathcal J_G^{(e)}}$ and $\sup_{J^{11}+J^{10}+J^{01}=J_G}$ involve   $(J^{11}+J^{01})p \choose J_e$ and $J_G+2 \choose 2$ terms, respectively, and we have   $$\log{J^{11}+J^{01})p \choose J_e}\leq J_e \log (e J_G p/J_e), \quad \log{J_G+2 \choose 2}\leq C\log(e J_G).$$
Applying the union bound gives the desired result. \qed

\subsection{Proof of Theorem \ref{th_minimax}}
The proof is based on the general reduction scheme  in Chapter 2 of \cite{tsybakov2009introduction}.
The key is to design proper {least favorable} signals  in different situations.
We consider two cases.

(i) $J_G\log (ep/J_G) \leq J_e\log (eJ_G^{2}/J_e) $.
Define a signal subclass
\begin{gather*}
\mathcal B^{1}(J_G, J_e) = \left\{\bsbOmega =[\bsbb, \bsbPhi]^\top: b_j = 1 \mbox{ if } j\in [J_G] \mbox{ and }  0 \mbox{ if } J_G< j\leq p,\right.\\
\left.\phi_{k,j} = 0 \mbox{ or }  \gamma R \mbox{ if } (k, j)\in [J_{G}]\times [J_G] \mbox{ and }  0 \mbox{ otherwise}, \bsbPhi=\bsbPhi^\top, J_e(\bsbPhi) \leq J_e \right\},
\end{gather*}
where $R= \frac{\sigma}{{\overline{\kappa}^{1/2}}  } (\log (e J_G^2 /J_e))^{1/2} $ and $\gamma>0$ is a small constant to be chosen later.
Because  $b_j=1$ ($1\leq j \leq J_G$),   any $\bsbOmega\in \mathcal B^{1}(J_G, J_e)$ satisfies SH, and thus ${\mathcal B}^{1}(J_G, J_e) \subset\text{{SH}}(J_G,J_e) $.
By Stirling's approximation, $$\log |\mathcal B^{1}(J_G, J_e)|= \log { J_G^2 \choose J_e}\geq   J_e \log (J_G^2/J_e)\geq c J_e \log ( eJ_G^2/J_e)$$ for some universal constant $c$.

Let $\rho(\bsbOmega_1, \bsbOmega_2)=\|\bsbOmega_1 - \bsbOmega_2\|_0$ be the Hamming distance.
 By Lemma 8.3 in   \cite{Rigollet11}, there exists  a subset ${\mathcal B}^{10}(J_G, J_e)\subset {\mathcal B}^{1}(J_G, J_e)$ such that
\begin{eqnarray*}
&\log | {\mathcal B}^{10}(J_G, J_e)| \geq c_1 J_e \log ( e J_G^2/J_e),
\\
&\rho(\bsbOmega_1, \bsbOmega_2) \geq c_2 J_e, \,\forall \bsbOmega_1, \bsbOmega_{2} \in \mathcal B^{10},  \bsbOmega_1\neq \bsbOmega_2
\end{eqnarray*}
for some universal constants $c_1, c_2>0$.
Then $\| \bsbOmega_1 - \bsbOmega_2\|_F^2 = \gamma^2 R^2 \rho(\bsbOmega_1, \bsbOmega_2) \geq c_2 \gamma^2 R^2 J_e $. It follows from the restricted conditional number assumption  that
\begin{align}
\|\breve \bsbX \vect(\bsbOmega_1) - \breve \bsbX \vect(\bsbOmega_2)  \|_2^2 \geq c_2 \underline{\kappa}   \gamma^2 R^2 J_e \label{separationLBound}
\end{align}
 for any $\bsbOmega_1, \bsbOmega_{2} \in \mathcal B^{10}$,  $\bsbOmega_1\neq \bsbOmega_2$, where $\underline{\kappa}/\overline{\kappa}$ is a positive constant.

For Gaussian models, the  Kullback-–Leibler divergence of $\mathcal N(\breve \bsbX \vect(\bsbOmega_2),\sigma^2 \bsb{I})$ (denoted by $P_{\bsbOmega_2}$) from $\mathcal N(\breve \bsbX \vect(\bsbOmega_1),\sigma^2 \bsb{I})$ (denoted by $\mathcal P_{\bsbOmega_1}$) is $\mathcal K(\mathcal P_{\bsbOmega_1}, \mathcal P_{\bsbOmega_2}) = \frac{1}{2\sigma^2} \|\breve\bsbX \vect(\bsbOmega_1) - \breve \bsbX \vect(\bsbOmega_2) \|_2^2$. Let $P_{\bsb{0}}$ be $\mathcal N(\bsb{0}, \sigma^2 \bsb{I})$. By the assumption again,
for any $\bsbOmega\in \mathcal B^1(J_G, J_e)$,
\begin{align*}
\mathcal K(P_{\bsb{0}}, P_{\bsbOmega}) \leq \frac{1}{2\sigma^2}\overline{\kappa}  \gamma^2 R^2 \rho(\bsb{0}, \bsbOmega) \leq \frac{\gamma^2}{\sigma^2}\overline{\kappa}  R^2 J_e,
\end{align*}
where we used $\rho(\bsbOmega_1, \bsbOmega_2)=J_e(\bsbPhi_1-\bsbPhi_2) \leq 2 J_e$. Therefore,
\begin{align}
\frac{1}{|\mathcal B^{10}|}\sum_{\bsbOmega\in \mathcal B^{10}} \mathcal K(P_{\bsb{0}}, P_{\bsbOmega})\leq \gamma^2 J_e\log (e J_G^2/J_e)\leq \frac{\gamma^2}{c} \log {J_G^2 \choose J_e}. \label{KLUBound}
\end{align}

Combining  \eqref{separationLBound} and \eqref{KLUBound} and choosing a sufficiently small value for  $\gamma$, we can apply Theorem 2.7 of \cite{tsybakov2009introduction} to get the desired  lower bound.

(ii) $J_G\log (ep/J_G) > J_e\log (eJ_G^2/J_e) $.
Define a signal subclass
\begin{gather*}
\mathcal B^{2}(J_G) = \{\bsbOmega =[\bsbb, \bsbPhi^\top]^\top: b_j \in \{0, 1\} \ \forall j \in  [ p], \bsbPhi= \bsb{0}, J_G(\bsbOmega)\leq J_G\},
\end{gather*}
where $R= \frac{\sigma}{{\overline{\kappa}^{1/2}}  } \{\log (e p/J_G)\}^{1/2} $ and $\gamma>0$ is a  small  constant. The afterward treatment is similar to   (i). The details are omitted.\\

For WH, we redefine  ${\mathcal B}^{1} $  as
$
\mathcal B^{1}(J_G, J_e) = \{\bsbOmega =[\bsbb, \bsbPhi^\top]^\top: b_j = 1$  if $j\in [J_G]$  and $ 0$  otherwise, $\phi_{k,j} =  0 \mbox{ or }  \gamma R$  if  $k\geq j$ and  $1\leq j \leq J_G$,   and $0$  otherwise,  $J_e^w(\bsbPhi) \leq J_e \}$,
where  $R= \frac{\sigma}{{\overline{\kappa}^{1/2}}  } (\log (e J_G p /J_e))^{1/2}$. Then ${\mathcal B}^{1}(J_G, J_e)\subset \mbox{WH}(J_G, J_e)$. The rest  follows the same lines as the proof for SH.

\section{Coordinatewise error bound for support recovery}
\label{sec:coord}
In this part, we show that GRESH estimators can recover the   sparsity pattern of the true signal, i.e.,   $\mathcal J^* =\{j:   \vect( {{\bm\Omega}^*})_j  \ne 0 \}$, with high probability. The result is  of the type of       \cite{lounici2008} and  \cite{ravikumar2010high},  but  we do not have to assume the stringent irrepresentable conditions   or mutual coherence conditions. For more discussions, see \cite{zhang2009some}.

\noindent {\sc Assumption $\mathcal{C}(\mathcal{J}_e,\mathcal{J}_G,\zeta, \lambda_1, \lambda_2)$}. Given $\mathcal{J}_e \subset [p^2]$,  $\mathcal{J}_G \subset [p]$,  $\zeta\ge 0$, $\lambda_1, \lambda_2>0$, for any $\bsbDelta = [\bsbDelta_{b},\bsbDelta_{\Phi}^{\top} ]^{\top} \in \mathbb{R}^{(p+1) \times p}$ satisfying  $\bsbDelta_{\Phi}=\bsbDelta_{\Phi}^{\top}$, the following inequality holds
  \begin{equation}
      \begin{aligned}
      &  \|\breve{\bsbX}\|_2(\|(\bsbDelta_{\Phi})_{\mathcal{J}_e}\|_1+\|(\bsbDelta_{\mathcal{J}_G})^{\top}\|_{2,1})  +   (J_e +J_G) \|\vect(\bm\Delta)\|_{\infty}^2 / (2\lambda_1\wedge 2\lambda_2)\\
      \le\ & \|\breve{\bsbX}\|_2 ( \|(\bsbDelta_{\Phi})_{\mathcal{J}_e^c}  \|_1+\|(\bsbDelta_{\mathcal{J}_G^c})^{\top}\|_{2,1})  +\zeta \sqrt{J_e +J_G}\|\bsbX\bsbDelta_{b}+\bar{\bsbX}\vect(\bsbDelta_{\Phi})\|_2.
      \end{aligned}\label{comparegcondInf}
      \end{equation}

      \begin{theorem}
  \label{th_infbound}
  Assume $\boldsymbol{\epsilon}\sim N(\bsb{0}, \sigma^2 \boldsymbol{I})$.   Take    $\lambda_1=A_1 \sigma \sqrt{\log (ep)}$ and $\lambda_2=A_2 \sigma \sqrt{\log (ep)}$  in   \eqref{typeA_crit_oracle}and suppose $(\bsbX, \bar{\bsbX},\breve{\bsbX})$ satisfies $\mathcal{C}(\mathcal{J}_e^*,\mathcal{J}_G^*,\zeta, \lambda_1, \lambda_2)$ for some $\zeta \ge 0$, where  $\mathcal J_e^* = \mathcal J_e(\bsbPhi^*)$, $\mathcal J_G^*=\mathcal J_G(\bsbb^*, \bsbPhi^*)$.   
Let $\hat{\bsbOmega}=\big[\hat{\bsb{b}},\hat{\bsbPhi}^{\top}\big]^{\top}$ be a global minimizer of \eqref{typeA_crit_oracle}. Define      $
      \hat{\mathcal J}  =  \{j: | \vect(\hat{\bm\Omega})_j| >      \sqrt{1+ \zeta^2} (\lambda_1 \wedge \lambda_2) \}
$,  and   $\mathcal J^* =\{j:   \vect( {{\bm\Omega}^*})_j  \ne 0 \}$.   Then for sufficiently large constants $A_1$ and $A_2$, with   probability at least $1-Cp^{-c (A_1^2\wedge A_2^2)}$,   the estimate  $\hat{\bm\Omega}$ satisfies
$$
      \|\vect(\hat{\bm\Omega}-\bm\Omega^*)\|_{\infty}  \le      \sqrt{1+ \zeta^2} (\lambda_1 \wedge \lambda_2).
 $$
 If, in addition, the minimum signal strength satisfies
 \begin{align}
      \min_{j\in\mathcal J^*  }| \vect(\bm\Omega^*)_j|   > 2    \sqrt{1+ \zeta^2} (\lambda_1 \wedge \lambda_2), \label{signalstrength}
\end{align}
       then   $\hat{\mathcal J} = \mathcal J^*$, with   probability at least $1-Cp^{-c (A_1^2\wedge A_2^2)}$.
\end{theorem}
To recover the true support, signal strength conditions like \eqref{signalstrength} must be imposed. Interestingly, compared with Theorem 5.1 in \cite{Loun2011}, hierarchical variable selection can accommodate smaller signals than group variable selection.
\begin{proof}

From the proof of   Theorem \ref{th_oracle}, we get the following inequality with probability at least $1- C p ^{-cL}$
  \begin{align*}
&\frac{1}{2} M(\bsbb^*-\hat{\bsbb},\bsbPhi^*-\hat{\bsbPhi})+A \lambda^o \rho\|(\bsbDelta_{\Phi})_{\mathcal{J}^c_e}\|_1+A\lambda^o\rho\|(\bsbDelta_{\mathcal{J}^c_G})^{\top}\|_{2,1}\\
\leq\ &\frac{1}{2} (1+\frac{3a_{1}}{A_0-6}) M(\bsbb^*-\bsbb,\bsbPhi^*-\bsbPhi)+\frac{1}{2} (\frac{1}{a_1}-1) M(\bsbDelta_{b},\bsbDelta_{\Phi}) \\
&+A \lambda^o \rho\|(\bsbDelta_{\Phi})_{\mathcal{J}_e}\|_1+ A \lambda^o \rho\|(\bsbDelta_{\mathcal{J}_G})^{\top}\|_{2,1} +\frac{a_1 L}{2(1-6/A_0)}\sigma^2 (J_e + J_G) \log (ep),
\end{align*}
where $\rho = \|\breve\bsbX\|_2$,  $L$ is a universal constant, $\lambda^o = \sigma\sqrt{\log (ep)}$,   $A^2 = A_0 L$, and  $a_1, A_0$ are positive constants to be determined. Also recall that
$\lambda_1 = A_1 \lambda^o$, $\lambda_2 = A_2 \lambda^o$, $A =  A_1  \wedge A_2 $.

By the regularity condition, we have
\begin{align*}
& \frac{1}{2}(J_e +J_G) \|\vect(\bm\Delta)\|_{\infty}^2 + \frac{1}{2} M(\bsbb^*-\hat{\bsbb},\bsbPhi^*-\hat{\bsbPhi})\\
\le\ & \frac{1}{2} (1+\frac{3a_{1}}{A_0-6}) M(\bsbb^*-\bsbb,\bsbPhi^*-\bsbPhi)+\frac{1}{2} (\frac{1}{a_1}-1) M(\bsbDelta_{b},\bsbDelta_{\Phi})\\&+ A \lambda^o \zeta  \sqrt{J_e +J_G}\|\bsbX\bsbDelta_{b}+\bar{\bsbX}\vect(\bsbDelta_{\Phi})\|_2 +\frac{a_1 L}{2(1-6/A_0)}\sigma^2 (J_e + J_G) \log (ep)\\
\le\ & \frac{1}{2} (1+\frac{3a_{1}}{A_0-6}) M(\bsbb^*-\bsbb,\bsbPhi^*-\bsbPhi)+\frac{1}{2} (\frac{1}{a_1}+\frac{1}{a_2}-1) M(\bsbDelta_{b},\bsbDelta_{\Phi})\\&+\left\{\frac{a_1 }{2(1-6/A_0)}+\frac{a_2 A_0 }{2} \zeta^2\right\}L\sigma^2 (J_e + J_G) \log (ep) ,
\end{align*}
for any $a_2>0$.
  Taking $\bm b = \bm b^*, \bm\Phi = \bm\Phi^*$ and choosing $a_1=a_2=1, A_0\ge 7$ give $  (J_e^* +J_G^*)\|\vect(\hat {\bm\Omega} - {\bm\Omega}^*)\|_{\infty}^2
  \le   A_{0} ({1 +  \zeta^2}) L\sigma^2 (J_e^* + J_G^*) \log (ep)  $. The conclusion follows.
\end{proof}

\section{An ADMM algorithm}
\label{appadmm}
We describe an ADMM algorithm for solving problem \eqref{general0}. Recall $\bar \bsbZ=[\bsbz_1 \odot \bsbz_1,  \cdots, \bsbz_1 \odot \bsbz_p, \cdots, \bsbz_p \odot \bsbz_1, \cdots, \bsbz_p \odot \bsbz_p]$,  $\breve \bsbZ=[\bsbx_1, \bsbz_1\odot\bsbz_1, \cdots, \bsbz_1\odot\bsbz_p, \cdots, \bsbx_p,\bsbz_p\odot\bsbz_1, \cdots, \bsbz_p\odot\bsbz_p]$, and $\bsb{\Lambda}_{\Omega} = \bsb{1}\bsb{\lambda}_{\Omega}^{\top}$. To apply ADMM, we rewrite \eqref{general0} into
\begin{equation}\label{general_ADMM}
\begin{split}
\min_{\bm{\Omega},\bm{\Gamma},\bm{\Upsilon}} &\frac{1}{2}\|\bm{y}-\bm{Xb}-\text{diag}(\bm{Z\Phi Z}^\top)\|_2^2 + \|\bm{\Lambda}_{\Omega}^{\top}\odot \bm{\Upsilon}^{\top}\|_{2,1}\\
& + \big\{\|\bm{\lambda}_b\odot\bm{\Gamma}_{b}\|_1 + \|\bm{\Lambda}_{\Phi}\odot \bm{\Gamma}_{\Phi}\|_1 \big\} \\
&\quad\quad\text{s.t. } \bm{\Gamma}_{\Phi}=\bm{\Gamma}_{\Phi}^{\top}, \bm{\Gamma} = \bm{\Omega}, \bm{\Upsilon} = \bm{\Omega},
\end{split}
\end{equation}
where  $\bm{\Upsilon} \in \mathbb{R}^{n\times (p^2+p)}$, $\bm{\Gamma} = [\bm{\Gamma}_b, \bm{\Gamma}_{\Phi}^{\top}]^{\top}$ with $\bm{\Gamma}_b \in \mathbb{R}^p, \bm{\Gamma}_{\Phi} \in \mathbb{R}^{p^2\times p^2}$.
The augmented Lagrangian of \eqref{general_ADMM} can be formed by use of two Lagrangian multiplier matrices ($\bm{L}_1,\bm{L}_2$) and a penalty parameter $\rho$:
\begin{align*}
\begin{split}
&L_{\rho}(\bm{\Omega},\bm{\Gamma},\bm{\Upsilon},\bm{L}_1,\bm{L}_2)\\
=~& \frac{1}{2}\|\bm{y} - \breve{\bm{Z}} \text{vec}(\bm{\Omega})\|_2^2 + \big\{\|\bm{\lambda}_b\odot\bm{\Gamma}_{b}\|_1 + \|\bm{\Lambda}_{\Phi}\odot \bm{\Gamma}_{\Phi}\|_1 + \iota_{\{\bm{\Gamma}_{\Phi}=\bm{\Gamma}_{\Phi}^{\top}\}} \big\} + \|\bm{\Lambda}_{\Omega}^{\top}\odot \bm{\Upsilon}^{\top}\|_{2,1}\\
~& + \rho\bm{1}^{\top}(\bm{L}_1\odot(\bm{\Omega}-\bm{\Gamma}) + \bm{L}_2\odot(\bm{\Omega}-\bm{\Upsilon}))\bm{1} + \frac{\rho}{2}(\|\bm{\Omega}-\bm{\Gamma}\|_F^2 + \|\bm{\Omega}-\bm{\Upsilon}\|_F^2).
\end{split}
\end{align*}

Based on the proof of Theorem \ref{th_conv},   we can solve the sub-optimization  problems of $\bm{\Gamma}$ and $\bm{\Upsilon}$  by some proximity operators; the details are omitted. The full ADMM algorithm is given as follows.

\begin{algorithm}
\caption{{\small An ADMM algorithm }\label{alg2}}
\begin{tabbing}
  \enspace \textbf{Inputs}: \\
  \enspace Data: $\bm{X}$, $\bm{Z}$, $\bm{y}$. Regularization parameters: $\bm{\lambda}_{b},\bm{\Lambda}_{\Phi},\bm{\lambda}_{\Omega}$.  \\
  \textbf{ repeat} \\
  \enspace \quad\quad 1. $\text{vec}(\bm{\Omega}) \leftarrow \big(\breve{\bm{Z}}^{\top}\breve{\bm{Z}} + 2\rho\bm{I}\big)^{-1}\big[\breve{\bm{Z}}^{\top}\bm{y} + \rho~\text{vec}(\bm{\Gamma}+\bm{\Upsilon} - \bm{L}_1-\bm{L}_2 )\big]$;\\
  \enspace \quad\quad 2. $\bm{\Upsilon}[:,k] \leftarrow \vec{\bm{\Theta}}_S(\bm{\Omega}[:,k] + \bm{L}_2[:,k]; \bm{\lambda}_{\Omega}[k]/\rho),~\forall k: 1\le k\le p $;\\
  \enspace \quad\quad 3. $\bm{\Omega}[2:\text{end},:] \leftarrow (\bm{\Omega}[2:\text{end},:] + \bm{\Omega}[2:\text{end},:]^{\top})/2$;\\
  \enspace \quad\quad 4. $\bm{\Gamma} \leftarrow \Theta_S(\bm{\Omega} + \bm{L}_1; [\bm{\lambda}_b, \bm{\Lambda}_{\Phi}^{\top}]^{\top}/\rho);$\\
  \enspace \quad\quad 5. $\bm{L}_2 \leftarrow \bm{L}_2 + \bm{\Omega} - \bm{\Upsilon}$;\\
  \enspace \quad\quad 6. $\bm{L}_1 \leftarrow \bm{L}_1 + \bm{\Omega} - \bm{\Gamma}$;\\
  \enspace \textbf{until} convergence \\
  \enspace \textbf{Output}: $\bm{\Omega}$
\end{tabbing}
\end{algorithm}

\bibliographystyle{Chicago}
\bibliography{greshbib}

\end{document}